\documentclass[english]{amsart}

\title{Companion forms over totally real fields}

\usepackage{amsmath} 
\usepackage{amssymb} 
\usepackage{amscd}
\usepackage{tabularx}
\usepackage[all,cmtip,line]{xy}

\newcommand{\F}{\mathbf{F}}
\DeclareMathOperator{\Frob}{Frob}

 \newcommand{\To}{\longrightarrow}
 \newcommand{\isoto}{\stackrel{\sim}{\To}}

 \newcommand{\M}{\mathcal{M}}

 \newcommand{\bigO}{\mathcal{O}}

 \newcommand{\R}{\mathbf{R}}
 \newcommand{\Z}{\mathbf{Z}}
 \newcommand{\A}{\mathbf{A}_\mathbf{Q}}
 \newcommand{\Q}{\mathbf{Q}}
 
 \newcommand{\p}{\mathfrak{p}}
  \newcommand{\n}{\mathfrak{n}}
\renewcommand{\l}{\mathfrak{l}}
  
 \newcommand{\C}{\mathbf{C}}

 \newcommand{\set}[1]{\left\{#1\right\}}

 \newcommand{\Gal}{\operatorname{Gal}}
 \newcommand{\GL}{\operatorname{GL}}
 \newcommand{\rk}{\operatorname{rk}}
 \newcommand{\Jac}{\operatorname{Jac}}
 \newcommand{\Res}{\operatorname{Res}}
 
\newcommand{\Spec}{\operatorname{Spec}}
\newcommand{\tr}{\operatorname{tr}}
\newcommand{\End}{\operatorname{End}}

\newcommand{\Lie}{\operatorname{Lie}}
\newcommand{\Pic}{\operatorname{Pic}}
\newcommand{\Ver}{\operatorname{Ver}}

\usepackage{amsthm}
\theoremstyle{plain}
\newtheorem*{thmn}{Theorem}
\newtheorem{ithm}{Theorem}
 \newtheorem{thm}{Theorem}[section]
 \newtheorem{corollary}[thm]{Corollary}
 \newtheorem{lemma}[thm]{Lemma}
 \newtheorem{prop}[thm]{Proposition}
 \theoremstyle{definition}
 \newtheorem{defn}[thm]{Definition}
 \theoremstyle{remark}
 
 \numberwithin{equation}{section}

\begin{document}
\author{Toby Gee} \maketitle

\begin{abstract}

We show that if $F$ is a totally real field in which $p$ splits
completely and $f$ is a mod $p$ Hilbert modular form with parallel
weight $2<k<p$, which is ordinary at all primes dividing $p$ and has
tamely ramified Galois representation at all primes dividing $p$,
then there is a ``companion form'' of parallel weight $k':=p+1-k$.
This work generalises results of Gross and Coleman--Voloch for
modular forms over $\Q$.
\end{abstract}

\section{Introduction}Theorems on ``companion forms'' were proved by Gross
(\cite{gro90}) under the assumption of some unchecked
compatabilities, and then reproved by Coleman and Voloch
(\cite{cv92}) without such assumptions. We generalise the methods of
Coleman and Voloch to totally real fields.

If $f\in S_k(\Gamma_1(N);\overline{\F}_p)$ is a mod $p$
cuspidal eigenform, where $p\nmid N$, there is a continuous, odd,
semisimple Galois representation
\[\rho_f:\Gal(\overline{\Q}/\Q)\To\GL_2(\overline{\F}_p)
\]
attached to $f$. A famous conjecture of Serre predicts that all
continuous odd irreducible mod $p$ representations should arise in
this fashion. Furthermore, the ``strong Serre conjecture'' predicts
a minimal weight $k_{\rho}$ and level $N_{\rho}$, in the sense that
$\rho\cong\rho_g$ for some eigenform $g$ of weight $k_{\rho}$ and
level $N_{\rho}$ (prime to $p$), and if $\rho\cong\rho_f$ for some
eigenform $f$ of weight $k$ and level $N$ prime to $p$ then
$N_{\rho}|N$ and $k\geq k_{\rho}$. Following the recent work of Khare, Wintenberger and Kisin, Serre's conjecture is now a theorem. However, very little is known about the natural generalisation of Serre's conjecture to totally real fields. A much earlier result (which is used in the work of Khare and Wintenberger) is the implication ``weak Serre'' $\implies$ ``strong Serre''. If one wishes to generalise Khare and Wintenberger's work to totally real fields, a first step would be to generalise this earlier work.

In solving the problem of weight optimisation it becomes necessary
to consider the companion forms problem; that is, the question of
when it can occur that we have $f=\sum a_nq^n$ of weight $2\leq
k\leq p$ with $a_p\neq 0$, and an eigenform $g=\sum b_n q^n$ of
weight $k'=p+1-k$ such that $na_n=n^k b_n$ for all $n$. Serre
conjectured that this can occur if and only if the representation
$\rho_f$ is tamely ramified above $p$. This conjecture has been
settled in most cases in the papers of Gross (\cite{gro90}) and
Coleman-Voloch (\cite{cv92}).

We generalise these results to the case of parallel weight Hilbert
modular forms over totally real fields $F$. We assume that $p$
splits completely in $F$. Our arguments are based
on those of \cite{cv92}, with several non-trivial and crucial
adjustments (see below). Subsequent to the writing of this paper, we found a different approach to these problems via deformation theory and modularity lifting theorems, which yields rather stronger results (see \cite{gee06}). However, while our main result is weaker than that of \cite{gee06}, we feel that the results contained herein are of independent interest, because of the light they shed on the results of \cite{cv92}, and because the techniques (for example, the construction of a $\theta$-operator on Hilbert modular forms on Shimura curves) may well be of more general use; see the remarks below.

Our main theorem is the following:

\begin{ithm}\label{thma}Let $F$ be a totally real field in which an odd prime $p$
splits completely. Let $\pi$ be a mod $p$ Hilbert modular form of
parallel weight $2<k<p$ and level $\n$, with $\n$ coprime to $p$.
Suppose that $\pi$ is ordinary at all primes $\p|p$, and that the
mod $p$ representation
$\overline{\rho}_\pi:\Gal(\overline{F}/F)\to\GL_2(\overline{\F}_p)$
is irreducible and is tamely ramified at all primes $\p|p$. Then
there is a companion form $\pi'$ of parallel weight $k'=p+1-k$ and
level $\n$ satisfying
$\overline{\rho}_{\pi'}\cong\overline{\rho}_\pi\otimes\chi^{k'-1}$, where $\chi$ is the $p$-adic cyclotomic character.
\end{ithm}

In the case $F=\Q$, the mod $p$ Galois representations associated
to modular forms may be found in the $p$-adic Tate modules of the
Jacobians of certain modular curves. For other totally real fields
one can often use the Jacquet-Langlands correspondence to
realise them in the Tate modules of the Jacobians of Shimura
curves associated to certain quaternion algebras (see
\cite{car862}). Although arguments on ``level lowering'' for mod
$p$ modular forms have been generalised to the totally real case
by arguing on these Shimura curves (cf. \cite{jar99}), it does not
seem to be possible to generalise the methods of \cite{gro90} or
\cite{cv92} to work on these curves, because they have no
``modular'' interpretation as PEL Shimura varieties. Indeed, it is
not even clear that there is a natural ``Hasse invariant'' on the
curves corresponding to quaternion algebras. We work instead on
certain unitary curves defined over degree two CM extensions $E$
of $F$, which are of PEL type. In order to base change the Hilbert
modular forms to automorphic forms on these Shimura curves it is
necessary to twist by certain grossencharacters of $E$; we then
obtain companion forms on the Shimura curves before twisting back
again.

In order to apply the ideas of \cite{cv92} it is necessary to have
semistable models of the Shimura curves we use, and we construct
these following \cite{km}. The arguments in \cite{cv92} depend
crucially on the use of $q$-expansions, which are not available to
us due to the lack of cusps on our Shimura curves. It has therefore
been necessary to construct arguments that apply more generally; we
do this by systematic use of expansions at supersingular points. We suspect that this technique will be of use in much more general situations than the one considered here, and it may well be possible to use these expansions to replace many other arguments in the literature which require the use of expansions at cusps. See, for example, section 4.3 for an example of these techniques. In addition, we make use of a $\theta$-operator on our Shimura curves (the operator $M$ of section 4.3); we expect this operator to prove as useful in this context as the classical one has proved to be on modular curves. For example, this operator should prove useful in proving results about the possible weights of non-ordinary Hilbert modular forms. These techniques may additionally prove useful in the study of Serre weight for higher-dimensional unitary groups, about which very little is known; such study would take place on higher-dimensional Shimura varieties which are defined in an analogous fashion to the curves considered here.

In section \ref{shimurasec} we construct the semistable models of
Shimura curves that we need. We work here in rather greater
generality than we require for the rest of the paper, assuming
nothing about the ramification of $p$ in $F$. Section \ref{auto}
contains the necessary background material on automorphic forms
and base change. In section \ref{construct} we construct Hecke
operators on Shimura curves and present our generalisations of the
arguments of \cite{cv92}. Finally, we present the proof of Theorem
\ref{thma} in section \ref{final}.

Where possible, our notation follows that of the original papers
\cite{car861} and \cite{cv92}. One possible point of confusion is
that the unitary groups we study are, loosely speaking, forms of
$\Res_{F/\Q}\GL_2\times\GL_1$. When we refer to ``the'' Galois
representation attached to such an automorphic form, we mean the
Galois representation attached to the $GL_2$-part, rather than some
twist of this by the character corresponding to the $GL_1$-part.
However, we will always ensure that the $GL_1$-part of our forms is
trivial at $\p$, so the local representation at $\p$ will in any
case be independent of any such twist.

This paper is almost entirely the author's PhD thesis under the
supervision of Kevin Buzzard, and it is a pleasure to thank him both
for suggesting the problem and for numerous helpful conversations.
The proof of the combinatorial result needed in Theorem \ref{121} is
based on an argument of Noam Elkies; any deficiencies in its
exposition are due to me. It is a pleasure to thank Fred Diamond,
Frazer Jarvis, and Richard Taylor for several helpful conversations.

\section{Shimura Curves}\label{shimurasec}
\subsection{Notation for Shimura Curves} Our notation follows that
of \cite{car861}. Let $F$ be a totally real field of degree $d>1$
over $\Q$, and denote by $\tau_1,\dots,\tau_d$ the infinite places
of $F$. Let $\p$ be a finite place of $F$, and let $\kappa$ be the
residue field of $F$ at $\p$, with cardinality $q$ and
characteristic $p$, and write $\bigO_{\p}$ for the ring of
integers of $F_{\p}$, the completion of $F$ at $\p$. Write also
$\bigO_{(\p)}$ for $F\cap\bigO_{\p}$ and $\bigO_{\p}^{nr}$ for the
completion of the ring of integers of $F_{\p}^{nr}$, the maximal
unramified extension of $F_{\p}$. Fix once and for all an
isomorphism $\C\isoto\overline{\Q}_p$; this isomorphism will be
used implicitly in our discussion of automorphic forms. Let the
primes of $F$ above $p$ be $\p_1,\p_2,\dots,\p_m$, where
$\p_1=\p$. Let $B$ be a quaternion algebra over $F$ which splits
at exactly one infinite place, say $\tau_1$, and suppose that $B$
splits at $\p$. We fix a maximal order $\mathcal{O}_B$ of $B$, and
choose an isomorphism between $\mathcal{O}_{B,v}$ and
$M_2(\mathcal{O}_v)$ at all finite places $v$ of $F$ where $B$
splits. We also choose an isomorphism between $B_{\tau_1}$ and $M_2(\R)$.

Define $G=\Res_{F/\Q}(B^\times)$, a reductive group over $\Q$.
Then if $K$ is a compact open subgroup of $G(\A^\infty)$ we define
the associated \emph{Shimura curve} to be
$M_K(\C)=G(\Q)\backslash\left(G(\A^\infty)\times(\C-\R)\right)/K$ (here $\A^\infty$ denotes the finite adeles).
By work of Shimura, $M_K(\C)$ has a canonical model $M_K$ over $F$
(see page 152 of \cite{car861}). Let $\Gamma$ be the restricted
direct product of the $(B\otimes F_v)^\times$ for all $v\neq\p$.
If $K=K_{\p}K^{\p}$, with $K_\p\subseteq\GL_2(\bigO_\p)$ and
$K^\p\subseteq\Gamma$, we write $M_{K_{\p},K^{\p}}$ for $M_K$, and
we extend this notation in an obvious fashion in the following paragraphs. Define
the groups

\[W_n(\p) = \left\{\begin{pmatrix}a&b\\c&d\end{pmatrix} \in \GL_2(\bigO_{\p}) \biggl\vert
\begin{pmatrix}a&b\\c&d\end{pmatrix}\equiv \begin{pmatrix}1&0\\0&1\end{pmatrix} (\operatorname{mod }\p^n)
 \right\},\] and write $M_{n,H}$ for $M_{W_n(\p),H}$. We also write $M_{0,H}$ for $M_{\GL_2(\bigO_{\p}),H}$.
Carayol (\cite{car861}) shows further that for $H$ sufficiently
small (depending on $i$) there are smooth models
$\mathbf{M}_{i,H}$ for the $M_{i,H}(\C)$ over $\mathcal{O}_{\p}$.

Define the group

\[\operatorname{bal}.U_1(\p) = \left\{\begin{pmatrix}a&b\\c&d\end{pmatrix} \in \GL_2(\bigO_{\p}) \biggl\vert
a-1\in\p,d-1\in\p,c\in\p\right\}.\]We construct integral models
for $\mathbf{M}_{\operatorname{bal}.U_1(\p),H}$, via constructing integral models
for the associated unitary curves. It is these unitary curves that
will be used in the main part of this paper.

Carayol also defines $T=\Res_{F/\Q}(\mathbf{G}_m)$, and for any
compact open subgroup $U\subseteq T(\A^\infty)$ the finite set
$\mathcal{M}_U(\C)=T(\Q)\backslash
\left(T(\A^\infty)\times\pi_0(T(\R))\right)/U$. This set has a natural
right action of $\Gal(\overline{\Q}/F)$, acting through
$\Gal(F^{ab}/F)$ via the inverse of the natural action of
$\pi_0(T(\Q)\backslash T(\A))$ and the isomorphism given by local
class field theory (normalised to take geometric Frobenius to a
uniformiser), which gives rise to a finite $F$-scheme together
with an action of $T(\A^\infty)$. Similarly, we
define $U_{\p}^n$ to be the subgroup of $\mathcal{O}_\p^\times$
consisting of units congruent to $1$ modulo $\p^n$, and for
$V\subseteq(\mathbf{A}_F^{\infty,\p})^*$ an open compact subgroup
we put $\mathcal{M}_{n,V}=\mathcal{M}_{U_\p^n \times V}$, and
define $\mathcal{M}_n=\varprojlim_V\mathcal{M}_{n,V}$. Then, for
example, $\mathcal{M}_0$ is isomorphic to $\Spec(F_{\p}^{nr})$.

\subsection{Carayol's Work}\label{carayol}The curves $M_K$ are not PEL Shimura curves, so in order to construct integral
models Carayol instead works with the Shimura curves associated to
certain unitary groups, and then uses results of Deligne to relate
these to the original curves.

Choose $\lambda<0$ in $\Q$ so that $K=\Q(\sqrt{\lambda})$ is split
at $p$, and define $E=F(\sqrt{\lambda})$. Fix a choice of square
root of $\lambda$ in $\C$, so that the embeddings
$\tau_i:F\hookrightarrow \R$ extend to embeddings
$\tau_i:E\hookrightarrow \C$; we always consider $E$ as a subfield
of $\C$ via $\tau_1$. Choose a square root $\mu$ of $\lambda$ in
$\Q_p$, so that the morphism $E\To F_p\oplus F_p$,
$x+y\sqrt{\lambda}\mapsto(x+y\mu,x-y\mu)$ extends to an
isomorphism
$$E\otimes \Q_p\isoto F_p\oplus
F_p\isoto(F_{\p_1}\oplus\dots\oplus
F_{\p_m})\oplus(F_{\p_1}\oplus\dots\oplus F_{\p_m})$$which gives
an inclusion of $E$ in $F_\p$ via
$$E\hookrightarrow E\otimes\Q_p\isoto F_p\otimes F_p\stackrel{\operatorname{pr}_1}{\To} F_p\stackrel{\operatorname{pr}_1}{\To}
F_\p.$$Let $z\To\bar{z}$ denote conjugation in $E$ with respect to
$F$. Put $D=B\otimes_F E$ and let $l\To\overline{l}$ be the
product of the canonical involution of $B$ with conjugation in
$E$. Choose $\delta\in D$ such that $\overline{\delta}=\delta$ and
define an involution on $D$ by
$l^*=\delta^{-1}\overline{l}\delta$. Choose $\alpha\in E$ such
that $\overline{\alpha}=-\alpha$. Then if $V$ denotes the
underlying $\Q$-vector space of $D$ with left action of $D$, we
have a symplectic form $\Psi$ on $V$ given by
$$\Psi(v,w)=\tr_{E/\Q}(\alpha\tr_{D/E}(v\delta w^*)).$$ Then $\Psi$
is a nondegenerate alternating form on $V$ satisfying, for all
$l\in D$,
$$\Psi(lv,w)=\Psi(v,l^*w).$$Let $G'$ be the reductive
algebraic group over $\Q$ such that for any $\Q$-algebra $R$,
$G'(R)$ is the group of $D$-linear symplectic similitudes of
$(V\otimes_\Q R,\Psi\otimes_\Q R)$; alternatively, we may define
$G'$ as follows, as in section 2.1 of \cite{car861}. Let
$T_E=\Res_{E/\Q}(\mathbf{G}_m)$. Let $\mathbf{S}$ denote
$\Res_{\C/\R}(\mathbf{G}_m)$, and let $U_E$ be the subgroup of
$T_E$ defined by the equation $z\overline{z}=1$. Then we can
define $G''=G\times_ZT_E$, and a morphism
$$G''=G\times_ZT_E\stackrel{\nu'}{\To}T\times U_E$$ by $(g,z)\mapsto(\nu(g)z\overline{z},z/\overline{z})$, where $\nu$ is the reduced norm on $G$. Then
if $T'$ is the subtorus $\mathbf{G}_m\times U_E$ of $T\times U_E$,
$G'$ is the inverse image under $\nu'$ of $T'$. This description
has the advantage of furnishing us with an isomorphism
$$G'(\Q_p)=\Q_p^*\times\GL_2(F_\p)\times(B\otimes_FF_{\p_2})^*\times\dots\times(B\otimes_FF_{\p_m})^*,$$
where the $\Q_p^*$ factor is given by $\nu(g)z\overline{z}$.

Carayol (\cite{car861}, section 2.1) defines a morphism
$h':\mathbf{S}\To G'_\R$, such that the $G'(\R)$-conjugacy class
of $h'$ may be identified with the complex upper half plane, and
the composition $\mathbf{S}\stackrel{h'}{\To}G'_\R\To\GL(V_\R)$
defines a Hodge structure of type $\set{(-1,0),(0,-1)}$ on $V_\R$.
We can (and do) choose $\delta$ so that $\Psi$ is a polarisation
for this Hodge structure.

Let $\bigO_D$ be a maximal order of $D$, corresponding to a
lattice $V_\Z$ in $V$. The above decomposition of $E\otimes \Q_p$
induces decompositions of $D\otimes\Q_p$ and $\bigO_D\otimes\Z_p$:
$$\begin{array}{ccccccccc}
\bigO_D\otimes\Z_p & =& \bigO_{D^1_1}&\oplus\dots\oplus &\bigO_{D^1_m}&\oplus &\bigO_{D^2_1}&\oplus\dots\oplus &\bigO_{D^2_m}\\
\bigcap &&\bigcap&&\bigcap&&\bigcap&&\bigcap\\
D\otimes\Q_p & =& {D^1_1}&\oplus\dots\oplus &{D^1_m}&\oplus
&{D^2_1}&\oplus\dots\oplus &{D^2_m}
\end{array}$$where each $D^k_j$ is an
$F_{\p_j}$-algebra isomorphic to $B\otimes_F F_{\p_j}$, and
$l\mapsto l^*$ interchanges $D^1_j$ and $D^2_j$. In particular
$D^1_1$ and $D^2_1$ are isomorphic to $\operatorname{M}_2(F_p)$.
We can, and do, choose $\bigO_D$, $\alpha$ and $\delta$ so that
the following conditions hold:
\begin{enumerate}\item $\bigO_D$ is stable under $l\mapsto l^*$.
\item Each $\bigO_{D^k_j}$ is a maximal order in $D^k_j$, and
$\bigO_{D^2_1}\hookrightarrow D^2_1=\operatorname{M}_2(F_\p)$
is identified with $\operatorname{M}_2(\bigO_\p)$.
\item $\Psi$ takes integer values on $V_\Z$.
\item $\Psi$ induces a perfect pairing $\Psi_p$ on
$V_{\Z_p}=V_\Z\otimes\Z_p$.
\end{enumerate}
Then every $\bigO_D\otimes\Z_p$-module $\Lambda$ admits a
decomposition as
$$\Lambda=\Lambda^1_1\oplus\dots\oplus\Lambda^1_m\oplus\Lambda^2_1\oplus\dots\oplus\Lambda^2_m$$with
$\Lambda^k_j$ an $\bigO_{D^k_j}$-module. The
$\bigO_{D^2_1}$-module $\Lambda^2_1$ decomposes further as the
direct sum of two $\bigO_p$-modules $\Lambda^{2,1}_1$ and
$\Lambda^{2,2}_1$, the kernels of the idempotents
$\bigl(\begin{smallmatrix}1&0\\0&0\end{smallmatrix}\bigr)$ and
$\bigl(\begin{smallmatrix}0&0\\0&1\end{smallmatrix}\bigr)$
respectively.

Let $X'$ be the $G'(\R)$-conjugacy class of $h'$; then for any
open compact subgroup $K'\subseteq G'(\A^\infty)$ we have a
Shimura curve over $\C$ $$M'_{K'}(\C)=G'(\Q)\backslash
\left(G'(\A^\infty)\times X'\right)/K'.$$By work of Shimura this
has a canonical smooth and proper model $M'_{K'}$ defined over
$E$, which represents the functor
$$\M_{K'}:\underline{E\mbox{-algebras}}\To\underline{\mbox{Sets}}$$where for any $E$-algebra
$R$, $\M_{K'}(R)$ is the set of isomorphism classes of quadruples
$(A,\iota,\lambda,\overline{\eta})$, such that
\begin{enumerate}\item $A$ is an abelian scheme of relative
dimension $4d$ over $R$, with an action $\iota:D\hookrightarrow
\End(A)\otimes\Q$ of $D$. This action induces an action of $E$ on
$\Lie(A)$, and for each $\tau_i$ we let
$\Lie_{\tau_i}(A)=\Lie(A)\otimes_{E,\tau_i}\C$. Then we require
that $\Lie_{\tau_i}(A)=0$ for all $i>1$.
\item $\lambda$ is a polarisation
of $A$ so that the Rosati involution sends $\iota(l)$ to
$\iota(l^*)$.
\item $\overline{\eta}$ is a class modulo $K'$ of symplectic $D$-linear
similitudes $\eta:\hat{V}(A)\To V\otimes\A^\infty$, where
$\hat{V}=\hat{T}\otimes \Q$ is the product of the Tate modules of
$A$ over all primes, with symplectic structure coming from the
Weil pairings.
\end{enumerate}

Our choice of embedding of $E$ into $F_\p$ allows us to base
change this model to $F_\p$, where we again denote it by
$M'_{K'}$. This again represents a moduli problem, and in fact
Carayol constructs an integral model for $M'_{K'}$ by describing a
moduli problem over $\bigO_\p$ which is represented by a smooth
and proper scheme $\mathbf{M}'_{K'}$, assuming that $K'$ is small
enough, and which satisfies $\mathbf{M}'_{K'}\otimes F_\p\isoto
M'_{K'}$. Using the above notation we now describe the moduli
problem represented by $\mathbf{M}'_{n,H'}$ for $H'\subseteq
\Gamma'$ sufficiently small, where
$$G'(\A^\infty)=\Q_p^*\times\GL_2(F_\p)\times\Gamma',$$ so that
$\Gamma'=G'(\A^{p,\infty})\times(B\otimes_F
F_{\p_2})^*\times\dots\times(B\otimes_F F_{\p_m})^*$. In fact,
$\mathbf{M}'_{n,H'}$ represents the functor
$$\boldsymbol{\M}'_{n,H'}:\underline{\bigO_\p\mbox{-algebras}}\To\underline{\mbox{Sets}}$$
where for any $\bigO_\p$-algebra $R$, $\boldsymbol{\M}'_{n,H'}(R)$
is the set of isomorphism classes of quintuples
$(A,\iota,\lambda,\eta_\p,\overline{\eta}^\p)$ such that
\begin{enumerate}\item $A$ is an abelian scheme of relative
dimension $4d$ over $R$, with an action
$\iota:\bigO_D\hookrightarrow \End_R(A)$ of $\bigO_D$ such that
\begin{enumerate}\item the projective $R$-module $\Lie^{2,1}_1(A)$ has
rank one, and $\bigO_\p$ acts on it via $\bigO_\p\hookrightarrow
R$,
\item for $j\geq 2$, $\Lie^2_j(A)=0$.
\end{enumerate}
\item $\lambda$ is a polarisation
of $A$ of degree prime to $p$ such that the Rosati involution
sends $\iota(l)$ to $\iota(l^*)$.
\item $\eta_\p$ is an isomorphism of $(\bigO_\p/\p^n)$-modules
$$\eta_\p:(A_{\p^n})^{2,1}_1\isoto(\p^{-n}/\bigO_\p)^2.$$
\item $\overline{\eta}^\p$ is a class of isomorphisms
$\eta^\p=\eta^\p_p\oplus\eta^p:T^\p_p(A)\oplus\hat{T}^p(A)\isoto
W^\p_p\oplus\hat{W}^\p$ modulo $H'$, with $\eta_p^\p$ linear and
$\eta^p$ symplectic, where $\hat{T}^p(A)$ is the product of the
Tate modules away from $p$,
$T^\p_p(A)=(T_p(A))^2_2\oplus\dots\oplus(T_p(A))^2_m$,
$\hat{W}^p=V_\Z\otimes\hat{\Z}^p$ and
$W_p^\p=(V_{\Z_p})^2_2\oplus\dots\oplus(V_{\Z_p})^2_m$.
\end{enumerate}

We have a short exact sequence $$1\To G_1\To
G'\stackrel{\nu'}{\To}T'\To 1$$where $G_1$ is the derived subgroup
of $G$, and thus also of $G'$. Then for any compact open subgroup
$U'\subseteq T(\A^\infty)$ we define the finite set
$$\mathcal{M}'_{U'}(\C)=T'(\Q)\backslash
\left(T'(\A^\infty)\times\pi_0(T'(\R))\right)/U'.$$ This set has a
natural right action of $\Gal(\overline{\Q}/E)$ defined via class
field theory as before, and we have a model over $E$, denoted
$\mathcal{M}'_{U'}$. For $K'\subseteq G'(\A^\infty)$ open and
compact this yields an $E$-morphism with geometrically connected
fibres
$$M'_{K'}\stackrel{\nu'}{\To}\mathcal{M}'_{\nu'(K')}.$$

For $H'$ sufficiently small so that $M'_{0,H'}$ exists, we see
that there is a universal abelian variety $A$ together with an
action of $\bigO_D$. Then for all $n$ we have a locally free
$(\bigO_\p/\p^n)$-scheme $E'_{n,H'}$ defined by
$$E'_{n,H'}=\left(A_{\p^n}\right)^{2,1}_1.$$ In section 3.3 of
\cite{car861} Carayol notes the equality
$$E'_{n,H'}=[M'_{n,H'}\times(\p^{-n}/\bigO_\p)^2]/\GL_2(\bigO_\p/\p^n);$$
we can also define a group $L'_n$ over $\mathcal{M}'_{0,V'}$ via
$$L'_n=[\mathcal{M}'_{n,V'}\times(\p^{-n}/\bigO_\p)]/(\bigO_\p/\p^n)^*.$$

Locally for the \'{e}tale topology on $M'_{0,H'}$ we see that
$E'_{n,H'}$ is isomorphic to the constant $(\bigO_\p/\p^n)$-module
$(\bigO_\p/\p^n)^2$, and $M'_{n,H'}$ is the scheme over
$M'_{0,H'}$ parameterising isomorphisms
$k_\p:E'_{n,H'}\cong(\p^{-n}/\bigO_\p)^2$. Upon choosing a
uniformiser $\p$ we have an isomorphism
$$\bigwedge_{\bigO_\p}^2(\p^{-n}/\bigO_\p)^2\stackrel{\det}{\isoto}(\p^{-2n}/\p^{-n})\stackrel{\p^n}{\isoto}(\p^{-n}/\bigO_\p)$$
which induces an isomorphism
$$\bigwedge_{\bigO_\p}^2 E'_n\isoto\nu'^*L'_n.$$ We denote the
resulting alternating $\bigO_\p$-bilinear pairing (the ``Weil
pairing'') by $$e'_n: E'_n\times E'_n\To\nu'^*L'_n.$$

These definitions are easily extended to $\mathbf{M}'_{0,H'}$.
Indeed, if $\mathbf{A}$ is the universal abelian scheme over
$\mathbf{M}'_{0,H'}$, we put
$$\mathbf{E}'_{n,H'}=\left(\mathbf{A}_{\p^n}\right)^{2,1}_1.$$ These fit together to give a divisible $\bigO_{\p}$-module $\mathbf{E}'_\infty$. Let
$\boldsymbol{\M}'_{n,V'}$ be the normalisation of $\M'_{n,V'}$ in
the ring of $\p$-integers of $E$. Then (\cite{car861}, section
8.3) there is a unique simultaneous extension of the groups
$L'_{n,V'}$ to finite locally free groups $\mathbf{L}'_{n,V'}$
over $\boldsymbol{\M}'_{n,V'}$ satisfying obvious compatablities,
and a unique extension of
$M'_{K'}\stackrel{\nu'}{\To}\mathcal{M}'_{\nu'(K')}$ to a morphism
$\mathbf{M}'_{K'}\stackrel{\boldsymbol{\nu}'}{\To}\boldsymbol{\mathcal{M}}'_{\nu'(K')}$.
In section 9.1 of \cite{car861} Carayol demonstrates that there is
a unique extension of $e'_n$ to an alternating $\bigO_\p$-bilinear
pairing
$$\mathbf{e}'_n:\mathbf{E}'_n\times_{\mathbf{M}'_{0,H'}}\mathbf{E}'_n\To\boldsymbol{\nu}'^*\mathbf{L}'_n.$$

At each geometric point $x$ of $\mathbf{M}'_{0,H'}\otimes\kappa$
we have a local divisible $\bigO_{\p}$-module
$\left.\mathbf{E}'_{\infty}\right\vert_x$. For every positive
integer $h$ there is a unique divisible formal $\bigO_{\p}$-module
$\Sigma_h$ over $\bar{\kappa}$ of height $h$ (see \cite{car861},
section 0.8). From the existence of an alternating pairing we
conclude that $\left.\mathbf{E}'_{\infty}\right\vert_x$ must be
self-dual, and up to isomorphism either

(1)$\left.\mathbf{E}'_{\infty}\right\vert_x\cong\Sigma_1\times(F_{\p}/\bigO_{\p})$,
or

(2)$\left.\mathbf{E}'_{\infty}\right\vert_x\cong\Sigma_2.$

As in the modular setting, we call $x$ \emph{ordinary} in the
first case and \emph{supersingular} in the second. From section
9.4 of \cite{car861} it follows that the set of supersingular
points is finite and nonempty.

For the purposes of our constructions (and in particular the
modular interpretation of Hecke operators) it is more natural to
work with the unitary curves; we will make use of results from
\cite{jar99}, which deals instead with the quaternionic curves,
but it is easy to see that the arguments in question carry over
unchanged from the quaternion algebra to the unitary setting.

\subsection{The $\operatorname{bal}.U_1(\p)$-Problem }
\begin{defn}Let $S$ be a scheme over $\mathbf{M}'_{0,H'}$, where $H'$ is sufficiently small for $\mathbf{E}'_{1,H'}$ to exist. For a subscheme $Q$ of $\left.\mathbf{E}'_1\right\vert_S$, let
$\mathcal{I}_Q$ be the ideal sheaf defining $Q$. If $P$ is a
section of $\left.\mathbf{E}'_1\right\vert_S$ and $\mathcal{K}$ a
locally free sub-$\bigO_{\p}/\p$-group scheme of
$\left.\mathbf{E}'_1\right\vert_S$, we say that $P$
\emph{generates} $\mathcal{K}$, and write $\mathcal{K}=\langle
P\rangle$, if $\mathcal{K}$ is the subscheme of
$\left.\mathbf{E}'_1\right\vert_S$ defined by the ideal
$\prod_{\lambda\in(\bigO_{\p}/\p)}\mathcal{I}_{\lambda P}$.
\end{defn}Note that a generator in our sense is an
``$(\bigO_\p/\p)$-generator'' in the sense of \cite{km}.
\begin{defn}A \emph{$\operatorname{bal}.U_1(\p)$-structure} on $S$, an
$\mathbf{M}'_{0,H'}$-scheme, is an f.p.p.f. short exact sequence
of $\bigO_{\p}$-group schemes on $S$

\[0\To \mathcal{K}\To\left.\mathbf{E}'_1\right\vert_S\To \mathcal{K'}\To0\]
such that $\mathcal{K}$, $\mathcal{K'}$ are both locally free of
rank $q$, together with $P\in \mathcal{K}(S)$, $P'\in
\mathcal{K'}(S)$ such that $\mathcal{K}=\langle P\rangle $,
$\mathcal{K'}=\langle P'\rangle $.\end{defn}

\begin{defn}Define the functor
\begin{align*}\boldsymbol{\mathcal{M}}'_{\operatorname{bal}.U_1(\p),H'} :
\underline{\mbox{Sch}/\mathbf{M}'_{0,H'}}&\To \underline{\mbox{Sets}},\\
T&\mapsto \operatorname{bal}.U_1(\p)(T):=\{\operatorname{bal}.U_1(\p)\text{-structures on }T
\}.\end{align*}
\end{defn}

\begin{defn}Define the functor
\[\mathcal{M}'_{\operatorname{\operatorname{bal}}.U_1(\p),H'} : \underline{\mbox{Sch}/{M'}_{0,H'}}\To
\underline{\mbox{Sets}},\]
\[T\mapsto\{P \text{ a nowhere-zero section of }\left.E'_1\right\vert_T, P'\text{ a nowhere-zero section of }\left.E'_1\right\vert_T/\langle P\rangle \}.\]
\end{defn}

\begin{lemma}The functors $\boldsymbol{\mathcal{M'}}_{\operatorname{\operatorname{bal}}.U_1(\p),H'}$ and
$\mathcal{M'}_{\operatorname{\operatorname{bal}}.U_1(\p),H'}$ agree on $M'_{0,H'}$-schemes.
\end{lemma}
\begin{proof} To give an $\boldsymbol{\mathcal{M}}'_{\operatorname{\operatorname{bal}}.U_1(\p),H'}$-structure is to give a section $P$ which generates a finite flat
sub-$\bigO_{\p}$-group scheme $\langle P\rangle$ of
$\left.\mathbf{E}'_1\right\vert_S$ of rank $q$, together with a
generating section $P'$ of
$\left.\mathbf{E}'_1\right\vert_S/\langle P\rangle$. If $S$ is a
scheme over $M'_{0,H'}$, $\left.\mathbf{E}'_1\right\vert_S =
\left.E'_1\right\vert_S$ which is \'{e}tale over $S$. Thus locally
in the \'{e}tale topology $\left.\mathbf{E}'_1\right\vert_S$ is
non-canonically isomorphic to $(\bigO_{\p}/\p)^2$, so that
$\langle P\rangle $, $\langle P'\rangle $ are both isomorphic to
$(\bigO_{\p}/\p)$, and thus $P$, $P'$ are both nowhere zero
sections killed by $\p$. The converse is clear.
\end{proof}

\begin{lemma}\label{etale}The $F$-scheme $M'_{\operatorname{\operatorname{bal}}.U_1(\p),H'}$ represents the functor $\mathcal{M}'_{\operatorname{\operatorname{bal}}.U_1(\p),H'}$.
\end{lemma}
\begin{proof}Let $\mathcal{M}'_{1,H'}$ be the functor:
\[\mathcal{M}'_{1,H'} : \underline{\mbox{Sch}/M'_{0,H'}}\To
\underline{\mbox{Sets}},\]
\[T\mapsto\{\text{pairs }(P,Q)\text{ of sections of
}\left.E'_1\right\vert_T\text{ over }T\text{ which trivialise
}\left.E'_1\right\vert_T\}.\] Then $\mathcal{M}'_{1,H'}$ is
represented by $M'_{1,H'}$. For any object $S$ of
$\mbox{Sch}/M'_{0,H'}$ we have an action of $\GL_2(\bigO_{\p}/\p)$ on
$\mathcal{M}'_{1,H'}(S)$ given by
\[(P,Q)\mapsto(P,Q)\left(\begin{array}{cc}a & b\\c &
d\end{array}\right)=(aP+cQ,bP+dQ).\] Furthermore this action is
clearly functorial, and the equivalence classes under the action
of the subgroup
\[\operatorname{\operatorname{bal}}.\widetilde{U}_1(\p) = \left\{\left(\begin{array}{cc}a & b\\c &
d\end{array}\right) \in \GL_2(\bigO_{\p}/\p) \biggl\vert
a-1\in\p,d-1\in\p,c\in\p\right\}\] give elements of
$\mathcal{M}'_{\operatorname{\operatorname{bal}}.U_1(\p),H'}(S)$ in the obvious way, so that we
have a map of moduli problems
$\mathcal{M}'_{1,H'}\To\mathcal{M}'_{\operatorname{\operatorname{bal}}.U_1(\p),H'}$.

Locally for the \'{e}tale topology we may complete any
$\operatorname{\operatorname{bal}}.U_1(\p)$-structure $(P,P')$ (note that the pair $(P,P')$
determines $\mathcal{K}$, $\mathcal{K'}$) to a pair
$(P,Q)\in\mathcal{M}_{1,H}(S)$ trivialising
$\left.E'_1\right\vert_S$, and so to a morphism $S\To M'_{1,H'}$
of $M'_{0,H'}$-schemes. Thus to give a $\operatorname{\operatorname{bal}}.U_1(\p)$-structure on
$S$ is to give a section $S\To M'_{1,H'}/\operatorname{\operatorname{bal}}.\widetilde{U}_1(\p) =
M'_{\operatorname{\operatorname{bal}}.U_1(\p),H'}$, and $\mathcal{M}'_{\operatorname{\operatorname{bal}}.U_1(\p),H'}$ is
represented by the $F$-scheme $M'_{\operatorname{\operatorname{bal}}.U_1(\p),H'}$.
\end{proof}

\begin{lemma}\label{407}The functor $\boldsymbol{\mathcal{M}}'_{\operatorname{\operatorname{bal}}.U_1(\p),H'}$ is represented by the scheme $\mathbf{M}'_{\operatorname{\operatorname{bal}}.U_1(\p),H'}$.
\end{lemma}
\begin{proof}
By Lemma 7.5 of \cite{jar99} the functor
$\boldsymbol{\mathcal{M}}'_{U_1(\p),H'}$ is representable, so it
suffices to prove that the moduli problem associating
$\operatorname{\operatorname{bal}}.U_1(\p)$-structures to $U_1(\p)$-structures is relatively
representable. The additional structure is the choice of a
generator for $\left.\mathbf{E}'_1\right\vert_S/\langle P\rangle$,
so the result follows by Proposition 1.10.13 of \cite{km}.
\end{proof}
\begin{lemma}\label{pt1}The scheme $\mathbf{M}'_{\operatorname{\operatorname{bal}}.U_1(\p),H'}$ is regular of
dimension two, and the projection map
$\mathbf{M}'_{\operatorname{\operatorname{bal}}.U_1(\p),H'}\To \mathbf{M}'_{0,H'}$ is finite and
flat.
\end{lemma}
\begin{proof}(cf. \cite{km} 5.1, \cite{jar99} Theorem 7.6). That the natural projection map
$\mathbf{M}'_{\operatorname{\operatorname{bal}}.U_1(\p),H'}~\to~\mathbf{M}'_{0,H'}$ is finite is
immediate from the finiteness of $\mathbf{E}'_1$, and the
observation that any point $P$ is a generator for only finitely
many $\mathcal{K}$. To see this, note that it suffices to prove
this over an algebraically closed base field; but then the result
follows from the proof of \cite{dri76}, Proposition 1.7, which
gives an explicit description of $\mathbf{E}'_1$, and from an
argument similar to the one on page 241 of \cite{dera}. We prove
flatness and regularity via the homogeneity principle of
\cite{km}. Let $U$ be the set of points $x$ of
$\mathbf{M}'_{0,H'}$ such that for any lift
$y\in\mathbf{M}'_{\operatorname{\operatorname{bal}}.U_1(\p),H'}$ of $x$ the local ring at $y$ is
regular and flat over the local ring at $x$. We prove that:

(H1) $U$ is open.

(H2) $M_{\operatorname{\operatorname{bal}}.U_1(\p),H'}$ is finite \'{e}tale over $M_{0,H'}$, so
in particular (as $M_{0,H'}$ is regular of dimension two) $U$
contains all of $M_{0,H'}$.

(H3) If $U$ contains an ordinary point of
$\mathbf{M}'_{0,H'}\otimes\kappa$ then it contains all ordinary
points of $\mathbf{M}'_{0,H'}\otimes\kappa$.

(H4) If $U$ contains a supersingular point of
$\mathbf{M}'_{0,H'}\otimes\kappa$  then it contains all
supersingular points of $\mathbf{M}'_{0,H'}\otimes\kappa$.

(H5) $U$ contains a supersingular point of
$\mathbf{M}'_{0,H'}\otimes\kappa$.

It will then be immediate that $U=\mathbf{M}'_{0,H'}$ (because $U$
is open and there are only finitely many supersingular points,
(H5) implies that $U$ also contains an ordinary point, and hence
all points by (H3) and (H4)), so the lemma will follow
(two-dimensionality being clear, as $\mathbf{M}'_{0,H'}$ is itself
regular two-dimensional).

The projection $\mathbf{M}'_{\operatorname{\operatorname{bal}}.U_1(\p),H'}\To
\mathbf{M}'_{0,H'}$ is finite so proper, so the complement of $U$
is closed, being the union of the two closed sets in
$\mathbf{M}'_{0,H'}$ which are the images of the closed subsets of
$\mathbf{M}'_{\operatorname{\operatorname{bal}}.U_1(\p),H'}$ at which
$\mathbf{M}'_{\operatorname{\operatorname{bal}}.U_1(\p),H'}$ is not regular or not flat over
$\mathbf{M}'_{0,H'}$. Thus (H1) is proved.

Since $M'_{1,H'}$ is certainly finite \'{e}tale over $M'_{0,H'}$,
so too is $M'_{\operatorname{\operatorname{bal}}.U_1(\p),H'}$, so (H2) is immediate.

The proofs of (H3) and (H4) are identical to those in \cite{jar99}
Theorem 7.6, but we include them for completeness' sake. Let $x$
be a closed point of $\mathbf{M}'_{0,H'}\otimes\kappa$. If $y$ is
a closed point of $\mathbf{M}'_{\operatorname{\operatorname{bal}}.U_1(\p),H'}$ lying above $x$,
then the map
$\bigO_{\mathbf{M}'_{0,H'},x}\To\bigO_{\mathbf{M}'_{\operatorname{\operatorname{bal}}.U_1(\p),H'},y}$
is flat if and only if the induced map
$\widehat{\bigO^{sh}}_{\mathbf{M}'_{0,H'},x}\To\widehat{\bigO^{sh}}_{\mathbf{M}'_{\operatorname{\operatorname{bal}}.U_1(\p),H'},y}$
is flat, and $\bigO_{\mathbf{M}'_{\operatorname{\operatorname{bal}}.U_1(\p),H'},y}$ is regular
if and only if
$\widehat{\bigO^{sh}}_{\mathbf{M}'_{\operatorname{\operatorname{bal}}.U_1(\p),H'},y}$ is regular
(\cite{egaiv}, IV 18.8.8 and 18.8.13). We are thus reduced, by the
argument on p.133 of \cite{km}, to considering the case where $x$
is a geometric point of the special fibre.

Let $(\widehat{\mathbf{M}'}_{0,H'})_{(x)}$ be the completion of
the strict henselisation of $\mathbf{M}'_{0,H'}$ at $x$. Then by
\cite{car861} 6.6
$\left.\mathbf{E}'_{\infty}\right\vert_{(\widehat{\mathbf{M}'}_{0,H'})_{(x)}}$
is the universal deformation of
$\left.\mathbf{E}'_{\infty}\right\vert_x$, so that the isomorphism
class of the map
$\widehat{\bigO^{sh}}_{\mathbf{M}'_{0,H'},x}\To\widehat{\bigO^{sh}}_{\mathbf{M}'_{\operatorname{\operatorname{bal}}.U_1(\p),H'},y}$
depends only on the universal deformation of
$\left.\mathbf{E}'_{\infty}\right\vert_x$, which in turn depends
only on whether $x$ is ordinary or supersingular, as required.

It remains to prove (H5). Accordingly, let $x$ be a supersingular
point of the special fibre of $\mathbf{M}'_{0,H'}$. Let
$\mathcal{C}$ denote the category of complete noetherian local
$\bigO_{\p}^{nr}$-algebras with residue field $\bar{\kappa}$. By
\cite{jar99} Theorem 4.5, the functor \begin{align*}\mathcal{C}
&\To \underline{\mbox{Sets}},\\
R&\mapsto\{\text{isomorphism classes of deformations of }
\left.\mathbf{E}'_{\infty}\right\vert_x\},\end{align*} is
represented by $\bigO_{\p}^{nr}[[t]]$.
\begin{lemma}The set
$\operatorname{\operatorname{bal}}.U_1(\p)(x)$ contains precisely one element.
\end{lemma}
\begin{proof}The set
$\operatorname{\operatorname{bal}}.U_1(\p)(x)$ is the set of triples $(\mathcal{K},P,P')$ with
$\mathcal{K}$ a finite flat $\bigO_{\p}$-subgroup scheme of
$\left.\mathbf{E}'_1\right\vert_{\bar{\kappa}}$,
$P\in\mathcal{K}(\bar{\kappa})$ a generator of $\mathcal{K}$, and
$P'$ a $\bar{\kappa}$-valued generator of
$\mathcal{K'}=\left.\mathbf{E}'_1\right\vert_x/\mathcal{K}$. But
$\left.\mathbf{E}'_1\right\vert_x$ is local (as we are in the
supersingular case), so $\mathcal{K}$ and $\mathcal{K'}$ are both
local, so $P=P'=0$. This gives us a unique
$\operatorname{\operatorname{bal}}.U_1(\p)$-structure, as required (again, by the proof of
Proposition 1.7 of \cite{dri76} $\Sigma_2$ has a unique sub
$(\bigO_\p/\p)$-module of rank $q$).
\end{proof}
Consequently, we see that the moduli problem of
$\operatorname{\operatorname{bal}}.U_1(\p)$-structures on
$\left.\mathbf{E}'_{\infty}\right\vert_{(\widehat{\mathbf{M}'}_{0,H'})_{(x)}}/\bigO_{\p}^{nr}[[t]]$
is represented by an affine scheme $\operatorname{Spec}(A)$ with
$A$ a local ring. By the argument of \cite{km} Proposition 5.2.2
we need only show that the maximal ideal of $A$ is generated by
two elements, which we can do by mimicking the proof of \cite{km}
Theorem 5.3.2. In fact, the required argument is formally
identical to that given in \cite{km}; again, we can work with
parameters $X(P)$, $X'(P')$ on the formal $\bigO_\p$-module, and
then apply Proposition 5.3.4 of \cite{km} in the case $p^n=q$.
\end{proof}

\subsection{Canonical Balanced Structures} Just as in the modular curve
case, the existence of an alternating form on $\mathbf{E}'_1$
allows us to define canonical $\operatorname{\operatorname{bal}}.U_1(\p)$-structures; these are the analogue of the moduli structures obtained by fixing a value of the Weil pairing in the elliptic curve case. Firstly,
we define
\begin{defn}A section $P$ of $\boldsymbol{\nu}'^*\mathbf{L}'_{1,H}|_S$, $S$ an
$\mathbf{M}'_{0,H'}$-scheme, is a \emph{generator} of
$\boldsymbol{\nu}'^*\mathbf{L}'_{1,H'}(S)$ if the subscheme of
$\boldsymbol{\nu}'^*\mathbf{L}'_{1,H'}(S)$ defined by the ideal
$\prod_{\lambda\in(\bigO_{\p}/\p)}\mathcal{I}_{\lambda P}$ is
$\boldsymbol{\nu}'^*\mathbf{L}'_{1,H'}(S)$.
\end{defn}

\begin{thm}The functor \begin{align*}\underline{\mbox{Sch}/\mathbf{M}'_{0,H'}}&\To
\underline{\mbox{Sets}},\\
S&\mapsto \{\text{generators of
}\boldsymbol{\nu}'^*\mathbf{L}'_{1,H'}(S) \}\end{align*} is
representable by an affine scheme ${\mathbf{L}'}_{1,H'}^\times$.
\end{thm}
\begin{proof}In fact, ${\mathbf{L}'}_{1,H'}^\times$ is obviously a closed
subscheme of the (finite, so) affine scheme
$\boldsymbol{\nu}'^*\mathbf{L}'_{1,H'}$.
\end{proof}

\begin{defn}Given a $\operatorname{\operatorname{bal}}.U_1(\p)$-structure $(\mathcal{K},P,P')$ we define
$\langle P,P'\rangle:=\mathbf{e}_1(P,Q)$ where $Q$ is any lift of
$P'$ to $\left.\mathbf{E}'_1\right\vert_S$ (note that this is well
defined because $\mathbf{e}_1(P,P)=1$; note also that $Q$ may in
fact be defined over a finite flat extension of $S$, but
$\mathbf{e}_1(P,Q)$ will still be defined over $S$ by descent).
\end{defn}
We can now define canonical $\operatorname{\operatorname{bal}}.U_1(\p)$-structures exactly as in
\cite{km} chapter 9:

\begin{defn}A \emph{$\operatorname{\operatorname{bal}}.U_1(\p)^{\text{can}}$-structure} on $S$, an
${\mathbf{L}'}_{1,H'}^\times$-scheme, is an f.p.p.f. short exact
sequence of $\bigO_{\p}$-group schemes on $S$

\[0\To \mathcal{K}\To\left.\mathbf{E}'_1\right\vert_S\To \mathcal{K'}\To0\]
such that $\mathcal{K}$, $\mathcal{K'}$ are both locally free of
rank $q$, together with $P\in \mathcal{K}(S)$, $P'\in
\mathcal{K'}(S)$ such that the
${\mathbf{L}'}_{1,H'}^\times$-structure on $S$ provided by $P$ is
the canonical one (i.e. the one that $S$ possesses as an
${\mathbf{L}'}_{1,H'}^\times$-scheme).
\end{defn}
The obvious modifications of the proofs of Proposition 9.1.7 and
Corollaries 9.1.8--9.1.10 of \cite{km} are also valid in our case.
Thus

\begin{thm} The functor \begin{align*}\underline{\mbox{Sch}/{\mathbf{L}'}_{1,H'}^\times}&\To
\underline{\mbox{Sets}},\\
S&\mapsto \{\operatorname{\operatorname{bal}}.U_1(\p)^{\text{can}}-\text{structures on }S
\}\end{align*}is represented by a scheme
${\mathbf{M}'}^\text{can}_{\operatorname{\operatorname{bal}}.U_1(\p),H'}$ regular and
equidimensional of dimension two.
\end{thm}
\begin{proof}Exactly as in Chapter 9 of \cite{km}.
\end{proof}

Note that after a change of base to the tamely ramified extension
of $F_\p$ determined by our choice of uniformiser $\p$, the scheme
$\boldsymbol{\nu}'^*\mathbf{L}'_{1,H'}$ is \'{e}tale,
corresponding in the usual fashion to $\p^{-1}/\bigO_\p$ together
with a Galois action. Indeed, from section 8.2 of \cite{car861} we
see that if $F_\p^0$ is the extension of $F_\p^{nr}$ corresponding
to $(\bigO_\p/\p)^*$ via class field theory, then
$\boldsymbol{\nu}'^*\mathbf{L}'_{1,H'}$ corresponds to
$\p^{-1}/\bigO_\p$, together with the action of
$\Gal(\overline{F}_\p/F_\p^{nr})$ given by the composition
$$\Gal(\overline{F}_\p/F_\p^{nr})\To\Gal(F_\p^0/F_\p^{nr})\cong(\bigO_\p/\p)^*.$$

\subsection{Igusa Curves}Our study of the special fibre of
$\mathbf{M}'_{\operatorname{\operatorname{bal}}.U_1(\p),H'}$ makes use of Igusa curves; these
are defined in just the same way as in the modular case. If $S$ is
a $\mathbf{M}'_{0,H'}\otimes\kappa$-scheme, then we have the
absolute Frobenius morphism $F_{\operatorname{abs}}:S\To S$, given
by $s\mapsto s^q$ on affine rings. We also in the usual way obtain
for any $S$-scheme $Z$ a relative Frobenius $F:Z\To Z^{(\sigma)}$,
where $Z^{(\sigma)}$ is the pullback of $Z$ via
$F_{\operatorname{abs}}:S\To S$. In particular, we have a map
$F:\left.\mathbf{E}'_1\right\vert_S\To\left.\mathbf{E}'_1\right\vert^{(\sigma)}_S$.

Let $(-)^D$ denote Cartier duality. Then we define the
\emph{Verschiebung}
$V:\left.\mathbf{E}'_1\right\vert^{(\sigma)}_S\To\left.\mathbf{E}'_1\right\vert_S$
to be the dual of
\[F:(\left.\mathbf{E}'_1\right\vert_S)^D\To(\left.\mathbf{E}'_1\right\vert^D_S)^{(\sigma)}=(\left.\mathbf{E}'_1\right\vert^{(\sigma)}_S)^D.\]

\begin{defn}Let $S$ be a $\mathbf{M}'_{0,H'}\otimes\kappa$-scheme. Then an \emph{Igusa
structure} on $S$ is a point
$P\in\left.\mathbf{E}'_1\right\vert^{(\sigma)}_S(S)$ which
generates the kernel of
$V:\left.\mathbf{E}'_1\right\vert^{(\sigma)}_S\To\left.\mathbf{E}'_1\right\vert_S$.
\end{defn}
\begin{lemma}The moduli problem
\begin{align*}\operatorname{Ig}:\underline{\mbox{Sch}/\mathbf{M}'_{0,H'}\otimes\kappa}&\To
\underline{\mbox{Sets}},\\
S&\mapsto\{\text{Igusa structures on S}\}\end{align*} is
representable by a finite flat
$\mathbf{M}'_{0,H'}\otimes\kappa$-scheme
$\mathbf{M}'_{\operatorname{Ig},H'}$, which is regular
one-dimensional, of rank $q-1$ over
$\mathbf{M}'_{0,H}\otimes\kappa$.
\end{lemma}
\begin{proof}Everything follows as in the proof of Theorem 6.1.1. of \cite{km} (see sections 8 and 9 of \cite{jar99}), except for regularity, which follows from an obvious
modification of the argument on page 363 of \cite{km}.
\end{proof}

\subsection{The Special Fibre} We now analyse the special fibre of
${\mathbf{M}'}^\text{can}_{\operatorname{\operatorname{bal}}.U_1(\p),H'}$ with the aid of the
crossing theorem of \cite{km}, 13.1.3. We will show that
${\mathbf{M}'}^\text{can}_{\operatorname{\operatorname{bal}}.U_1(\p),H'}\otimes\kappa$ is the
union of two smooth curves crossing transversally above the
supersingular points of $\mathbf{M}'_{0,H'}\otimes\kappa$. As in
the modular case, one of these curves is essentially a copy of
$\mathbf{M}'_{0,H'}\otimes\kappa$, and corresponds to the
generators of $\ker F$, and the other is an Igusa curve, the
scheme of generators of $\ker V$.

There is a unique morphism ${\mathbf{L}'}_{1,H'}^\times\To\kappa$
(this can be seen explicitly by consideration of the construction
of tame extensions via adjoining roots of uniformisers) so the
reduction mod $\p$ of the $\operatorname{\operatorname{bal}}.U_1(\p)^{can}$ problem is given by
\begin{align*}\underline{\mbox{Sch}/\mathbf{M}'_{0,H'}\otimes\kappa}&\To \underline{\mbox{Sets}}\\
S&\mapsto\{\operatorname{\operatorname{bal}}.U_1(\p)\text{-structures }(\mathcal{K},P,P')\text{
with }\langle P,P'\rangle=1\}.\end{align*}The condition that the
pairing be $1$ is a closed condition, so this problem is
representable, say by $\mathbf{M}'_{\operatorname{\operatorname{bal}}.U_1(\p);\det=1,H'}$.

We need to analyse the possible $(\operatorname{\operatorname{bal}}.U_1(\p);\det=1)$-structures
on geometric points of $\mathbf{M}'_{0,H'}\otimes\kappa$. From
section 9 of \cite{jar99} we see that if $x$ is supersingular the
only possibility for $\mathcal{K}$ is $\ker F$, and if $x$ is
ordinary then $\mathcal{K}$ is either $\ker F$ or $\ker V$.
Furthermore:

\begin{lemma}Let $S$ be a $\mathbf{M}'_{0,H'}\otimes\kappa$-scheme. Let $P$ be a generator of $\ker V$. Then the triple $(\ker V,P,0)$ is a $(\operatorname{\operatorname{bal}}.U_1(\p);\det=1)$-structure on $S^{(\sigma)}$
and the triple $(\ker F,0,P \operatorname{mod} \ker F)$ is a
$(\operatorname{\operatorname{bal}}.U_1(\p);\det=1)$-structure on $S$. Furthermore, these
constructions define closed immersions
$\mathbf{M}'_{\operatorname{Ig},H'}\hookrightarrow
\mathbf{M}'^{(\sigma)}_{\operatorname{\operatorname{bal}}.U_1(\p);\det=1,H'}$ and
$\mathbf{M}'_{\operatorname{Ig},H'}\hookrightarrow
\mathbf{M}'_{\operatorname{\operatorname{bal}}.U_1(\p);\det=1,H'}$ respectively.
\end{lemma}
\begin{proof}That these define $\operatorname{\operatorname{bal}}.U_1(\p)$-structures on $S$ follows from the proof of \cite{jar99} Theorem 10.2. We need to check that
$\langle P,0\rangle=1$; but by definition $\langle
P,0\rangle=\mathbf{e}_1(P,P)=1$, as required.

Our claimed immersions are certainly
$\mathbf{M}'_{0,H'}\otimes\kappa$-maps between finite
$\mathbf{M}'_{0,H'}\otimes\kappa$-schemes, so are finite and hence
proper. We claim that they are injective on $S$-valued points for
all $\mathbf{M}'_{0,H'}\otimes\kappa$-schemes; this is immediate,
since we can clearly recover $P$ from the image. But (\cite{egaiv}
18.12.6) a proper monomorphism is a closed immersion, as required.
\end{proof}

\begin{thm}\label{pt2}$\mathbf{M}'_{\operatorname{\operatorname{bal}}.U_1(\p);\det=1,H'}$ is the disjoint union with transverse crossings at the supersingular points of two
smooth $\kappa$-curves.
\end{thm}
\begin{proof}This follows from an application of the Crossings
Theorem (\cite{km}, Theorem 13.1.3) and the previous lemma; one of
our curves is $\mathbf{M}'_{\operatorname{Ig},H'}$ and the other
is $\mathbf{M}'^{(\sigma^{-1})}_{\operatorname{Ig},H'}$. The
required argument is very similar to that in section 9 of
\cite{jar99}. That the crossings occur at the supersingular points
is simply the statement that $\ker F = \ker V$ precisely in the
supersingular case.

We claim that the completion of the strict henselisation of the
local ring at the supersingular points is
$\bar{\kappa}[[x,y]]/(xy)$, where as before $x=X(P)$, $y=X'(P')$,
from which it is immediate that the crossings at the singular
points are transverse.

From the crossings theorem it suffices to show that our two
components are given by $x=0$ and $y=0$ respectively. But one
component is given by $P=0$, which occurs if and only if
$x=X(P)=0$, and similarly the other component is defined by $y=0$,
as required.
\end{proof}
In fact, we claim that the universal formal deformation at a
supersingular point is given by $F_\p^0[[x,y]]/(xy-\pi)$, where
$\pi$ is a uniformiser for $F_\p^0$. The proof of this again
follows exactly as in the proof of theorem 5.3.2 of \cite{km},
which shows that the universal formal deformation is given by
$F_\p^0[[x,y]]/(f)$ for some $f$, and that the maximal ideal of
$F_\p^0[[x,y]]$ is generated by $x$, $y$ and $f$; the additional
information that the strict henselisation of the local ring mod
$p$ is $\bar{\kappa}[[x,y]]/(xy)$ gives the required result.

\section{Automorphic Forms}\label{auto}
\subsection{Definitions} We say that $\pi$ a cuspidal automorphic representation of $\GL_2/F$ is a
Hilbert modular form of weight $k\geq 2$ if for each
$\tau:F\hookrightarrow\R$ the representation $\pi_\tau$ is the
$(k-1)$-st lowest discrete series representation of $\GL_2(\R)$ with
central character $a\mapsto a^{2-k}$ (note that we are only working
with parallel weight modular forms). We will also say that
$\pi_\infty$ is of weight $k$ when this condition holds.

\begin{defn}A \emph{mod $p$ Hilbert modular form} is an
equivalence class of Hilbert modular forms, where two forms $\pi$,
$\pi'$ are equivalent if they have the same mod $p$ Galois
representations (see below); equivalently, if for all finite places
$v\nmid p$ of $F$ at which $\pi_v$ and $\pi'_{v}$ are unramified
principal series, say $\pi_v=\pi(\psi_1,\psi_2)$,
$\pi'_v=\pi(\psi'_1,\psi'_2)$, we have an equality
$$\{\psi_1(\Frob_v),\psi_2(\Frob_v)\}=\{\psi'_1(\Frob_v),\psi'_2(\Frob_v)\}\text{
mod }p.$$
\end{defn}

Let $\mathfrak{n}$ be an ideal of $\bigO_F$, and let
$U_1(\mathfrak{n})$ denote the subgroup of
$\prod\GL_2(\bigO_{F,v})$ consisting of elements
$\bigl(\begin{smallmatrix}a&b\\c&d\end{smallmatrix}\bigr)$ with
$c\in\mathfrak{n}$ and $(a-1)\in\mathfrak{n}$. Then (by the theory
of the conductor) for any Hilbert modular form $\pi$ there exists
some $\mathfrak{n}$ for which $\pi^{U_1(\mathfrak{n})}\neq 0$, and
we say that $\pi$ has level $\mathfrak{n}$. Note that a modular
form has infinitely many levels, all divisible by the minimal
level, the product of the local conductors of $\pi$. We say that a
mod $p$ form has weight $k$ and level $\mathfrak{n}$ if some form
in the equivalence class defining it does.

We have, as usual, a notion of a Hilbert modular form being ordinary
at $\p_i$; if  $\pi_{\p_i}$ is unramified with Satake parameters
$\alpha$, $\beta$, then we let $a_{\p_i}=q^{1/2}(\alpha+\beta)$.
Then $a_{\p_i}$ is the eigenvalue at $\p_i$ of the classical Hilbert
modular form corresponding to $\pi$, and we say that $\pi$ is
ordinary if $a_{\p_i}\neq 0$ mod $p$ for all $\p_i|p$, where $a_{\p_i}$ is the eigenvalue of $T_{\p_i}$. There is an
obvious notion of an ordinary mod $p$ form; we simply demand that it
is the reduction mod $p$ of a form which is ordinary at $\p$ for all $\p|p$. Then:

\begin{lemma}An ordinary mod $p$ Hilbert modular form of level
$\mathfrak{n}$ prime to $p$ is an ordinary mod $p$ form of weight
$2$ and level $\mathfrak{n}p$.
\end{lemma}
\begin{proof}This is an easy application of Hida theory, and in fact follows at once from Theorem 3 of
\cite{wil88}.
\end{proof}

Let $\pi$ be a Hilbert modular form, and let $M$ be the field of
definition of $\pi$ (that is, the fixed field of the automorphisms
$\sigma$ of $\C$ satisfying $\sigma\pi\cong\pi$). It is known that
$M$ is either a totally real or a CM number field, and that for
each prime $p$ of $\Q$ and embedding
$\lambda:M\hookrightarrow\overline{\Q}_p$ there is a continuous
irreducible representation
$$\rho_{\pi,\lambda}:\Gal(\overline{F}/F)\to\GL_2(M_\lambda)$$determined
(thanks to the Cebotarev density theorem and the Brauer-Nesbitt
theorem) by the following property: if $v\nmid p$ is a place of $F$
such that $\pi_v$ is unramified then
$\rho_{\pi,\lambda}|_{\Gal(\overline{F}_v/F_v)}$ is unramified with
$\Frob_v$ having minimal polynomial
$$X^2-t_vX+(\mathbf{N}v)s_v$$where $t_v$ is the eigenvalue of the
Hecke operator
$$\left[\GL_2(\bigO_{F_v})\begin{pmatrix}\varpi_v&0\\0&1\end{pmatrix}\GL_2(\bigO_{F_v})\right]$$
(where $\varpi_v$ is a uniformiser of $\bigO_{F_v}$) on
$\pi^{\GL_2(\bigO_{F_v})}$ and $s_v$ is the eigenvalue of the
operator
$$\left[\GL_2(\bigO_{F_v})\begin{pmatrix}\varpi_v&0\\0&\varpi_v\end{pmatrix}\GL_2(\bigO_{F_v})\right].$$As
usual, by the compactness of $\Gal(\overline{F}/F)$ we may
conjugate $\rho_{\pi,\lambda}$ to a representation valued in
$\GL_2(\bigO_{M,\lambda})$, and then reduce modulo the maximal
ideal of $\bigO_{M,\lambda}$ to get a continuous representation to
$GL_2(\overline{\mathbf{F}}_p)$. The semisimplification of this
representation is well-defined, and we denote it by
$\overline{\rho}_{\pi,p}$.

In particular, the above discussion shows that there is a
continuous mod $p$ Galois representation $\overline{\rho}_\pi$
canonically associated to any mod $p$ Hilbert modular form $\pi$,
determined as above by the characteristic polynomials of
Frobenius.

Let $\pi$ be a weight $2$ level $\mathfrak{n}p$ Hilbert modular
form, corresponding to our initial weight $k$ level $\mathfrak{n}$
mod $p$ Hilbert modular form. In order to apply the level lowering
results below, we need to assume that it is not the case that
$[F(\zeta_p):F]=2$ and
$\overline{\rho}_\pi|_{\Gal(\overline{F}/F(\zeta_p))}$ is
reducible; in the case that this does not hold it is easy to
construct the companion form ``by hand'', so from now on we ignore
this case until we treat it in the proof of Theorem A. At various
places in our arguments on Shimura curves we will need to assume
that the level is sufficiently large (equivalently, the compact
open subgroup of $G'(\mathbf{A}_\Q^\infty)$ corresponding to the
level structure is sufficiently small). This may be accomplished
via a trick originally due to Diamond and Taylor; namely, we
choose a prime $\mathfrak{q}\nmid\mathfrak{n}p$ such that there
are no congruences between forms of level $U_1(\mathfrak{n})$ and
$\mathfrak{q}$-new forms of level dividing $U_1(\mathfrak{n})\cap
U_1^1(\mathfrak{q})$, where $U_1^1(\mathfrak{q})$ is the subgroup
of $\prod\GL_2(\bigO_{F,v})$ consisting of elements
$\bigl(\begin{smallmatrix}a&b\\c&d\end{smallmatrix}\bigr)$ with
$c,a-1,d-1\in\mathfrak{q}$. We then work throughout with an
auxillary $U^1_1(\mathfrak{q})$-level structure, which we can
remove at the end due to the lack of congruences. There are
infinitely many such  primes $\mathfrak{q}$; see the remark
following Lemma 12.2 of \cite{jar99} for a proof of this.

If $d$ is even, we assume that there is a finite place $z\nmid p$ of
$F$ such that $\pi_z$ is not principal series. Fortunately, this
does not entail a loss of generality. By Theorem 1 of \cite{tay89}
there is a finite place $z\nmid\mathfrak{n}p$ where
$\overline{\rho}_\pi$ is unramified,
$\mathbf{N}_{F/\Q}(z)\equiv-1\text{ mod }p$, and a Hilbert modular
form $\tilde{\pi}$ of the same weight as $\pi$ and level
$U_1(\mathfrak{n})\cap U_z$ (where $U_z$ is the subgroup of
$\prod\GL_2(\bigO_{F,v})$ consisting of elements
$\bigl(\begin{smallmatrix}a&b\\c&d\end{smallmatrix}\bigr)$ with
$z|c$) such that $\tilde{\pi}_z$ is unramified special and
$\overline{\rho}_{\tilde{\pi}}\cong\overline{\rho}_\pi$. We can then
work throughout with $\tilde{\pi}$ in place of $\pi$, and remove the
auxillary level $z$ structure at the end thanks to Theorem A of
\cite{fuj99}.

We now show how to construct a holomorphic differential
corresponding to $\pi$ on a unitary Shimura curve. We choose a
quadratic imaginary extension $K/\Q$ as in section \ref{carayol},
but in addition to requiring it to split $p$, we require that it
splits all primes $l|N_{F/\Q}\mathfrak{n}$, and that it is disjoint from the
extension of $F$ given by the kernel of the Galois representation
$\overline{\rho}_\pi$ (this last condition ensures that the base
change of $\pi$ to $E=KF$ is cuspidal; see Theorem
\ref{cyclicbasechange} below). We also require that all ramified
places of $F$ and $B$ split in $E$. There is considerable freedom in the choice of $K$, and we will later exploit this to show that our final construction of a companion form is independent of any choices.

\begin{defn}If $\pi$ is an automorphic cuspidal representation of
$\GL_2(\mathbf{A}_F)$, and $\Pi$ is an automorphic representation
of $\GL_2(\mathbf{A}_E)$, then $\Pi$ is a \emph{base change lift}
of $\pi$, denoted $BC_{E/F}(\pi)$, if for every place $v$ of $F$,
and every $w|v$, the Langlands parameter attached to $\Pi_w$
equals the restriction to $W_{E_w}$ of the Langlands parameter
$\sigma_v:W_{F_v}\to\GL_2(\C)$ of $\pi_v$. In particular, if $v$
splits in $E$ then $\Pi_w\cong\pi_v$.
\end{defn}
\begin{thm}\label{cyclicbasechange}\begin{enumerate}\item Every cuspidal representation $\pi$ of
$\GL_2(\mathbf{A}_F)$ has a unique base change lift to
$\GL_2(\mathbf{A}_E)$; the lift is itself cuspidal unless $\pi$ is
monomial  of the form
$\pi\left(\operatorname{Ind}_{W_E}^{W_F}\theta\right)$. \item If
$\pi$, $\pi'$ have the same base change lift to
$\GL_2(\mathbf{A}_E)$ then $\pi'\cong\pi\otimes\omega$ for some
character $\omega$ of $F^\times
N_{E/F}(\mathbf{A}_E^\times)\backslash\mathbf{A}_F^\times$. \item
A cuspidal representation $\Pi$ of $\GL_2(\mathbf{A}_E)$ equals
$BC_{E/F}(\pi)$ for some $\pi$ if and only if $\Pi$ is invariant
under the action of $\Gal(E/F)$.\end{enumerate}
\end{thm}
\begin{proof}\cite{lan80}.
\end{proof}

Let $\pi_E=BC_{E/F}(\pi)$. Let $B$ denote the indefinite
quaternion algebra over $F$ which if $d$ is odd is ramified only
at $\tau_2,\dots,\tau_d$, and which if $d$ is even is ramified
precisely at $\tau_2,\dots,\tau_d$ and some $z\nmid\mathfrak{n}p$
with $\pi_z$ not principal series.

We now, as in \cite{br89}, twist $\pi_E$ by a character $\eta$ so
that $\eta\otimes(\eta\circ c)=\chi_\pi^{-1}\circ N_{E/F}$, where
$c$ denotes the nontrivial element of $\Gal(E/F)$ (the existence of
such characters is easily deduced from the arguments in the proof of
Lemma VI.2.10 of \cite{ht01}). Then if $^\vee$ denotes the
contragredient representation we have
\begin{align*}(\pi_E\otimes\eta)^\vee\circ c&\cong(\pi_E\circ
c)^\vee\otimes(\eta\circ c)^\vee\\&\cong\pi_E^\vee\otimes(\eta\circ
c)^\vee\\&\cong(\pi_E\otimes\chi_{\pi_E}^{-1})\otimes(\eta\circ
c)^{-1}\\&\cong\pi_E\otimes(\chi_\pi^{-1}\circ
N_{E/F})\otimes(\eta\circ
c)^{-1}\\&\cong\pi_E\otimes\eta.\end{align*} Now,
$(\chi_\pi)_\infty=1$, so we may suppose that $\eta$ is trivial at
the archimedean places. Furthermore, we choose $\eta$ so that
$\eta_{\p_i}$ is trivial at all $\p_i|p$, where we identify each
$\p_i$ with the place of $E$ given by our fixed choice of a place of
$K$ dividing $p$.

Then by the Jacquet-Langlands theorem (Theorem \ref{jl}) there is
a unique irreducible automorphic representation $\Pi$ of
$(D\otimes \A)^\times$ (where $D=B\otimes E$) such that
$JL(\Pi)=\pi_E\otimes\eta$. Furthermore (cf. page 199 of
\cite{ht01}) the condition that $(\pi_E\otimes\eta)^\vee\circ
c\cong\pi_E\otimes\eta$ gives $\Pi^*\cong\Pi$.

Now from Theorem VI.2.9 and Lemma VI.2.10 of \cite{ht01} we see that
there is a character $\psi$ of $\mathbf{A}^\times_K/K^\times$ and an
automorphic representation $\pi_1$ of $G'(\mathbf{A}_\Q)$ such that
$BC(\pi_1)=(\psi,\Pi)$ and $(\psi)_\infty$ is trivial. Furthermore,
from Theorem \ref{172} below we see that $(\pi_1)^{K'}\neq0$, where
$K'$ is defined as follows: $K'_l=G'(\Z_l)$ for all
$l\nmid\mathfrak{n}p$. For $l|\mathfrak{n}$ we have an isomorphism
$G'(\Q_l)=\Q_l^*\times\prod_{\l|l}\GL_2(F_\l)$, and we put
$K_l'=\Z_l\times\prod_{\l|l}\left\{\bigl(\begin{smallmatrix}a&b\\c&d\end{smallmatrix}\bigr)\in\GL_2(\bigO_\l):
c,d-1\equiv 0\text{ mod }\l^{v_\l(\mathfrak{n})}\right\}$. Finally,
$K_p'=\Z_p\times\prod_{\p_i|p}\operatorname{\operatorname{bal}}.U_1(\p_i)$. We will sometimes
denote the away-from-$\p$ part of the level by $H'$.

\subsection{Base Change}We now recall various results on
automorphic forms due to, amongst others, Jacquet, Langlands,
Clozel and Kottwitz. A convenient reference for these results is
\cite{ht01}, where they are all stated in far greater generality;
the reader should note, however, that many of the results quoted
have simpler proofs in our special cases than those in
\cite{ht01}.

\begin{thm}\label{jl}(Jacquet-Langlands) If $\rho$ is an irreducible
automorphic representation of $(D\otimes\A)^\times$ then there is
a unique automorphic representation $JL(\rho)$ of
$GL_2(\mathbf{A}_E)$ which occurs in the discrete spectrum and for
which $JL(\rho)^{S(D)}\cong \rho^{S(D)}$, where $S(D)$ is the set
of places at which $D$ ramifies. The image of $JL$ is the set of
irreducible automorphic representations $\pi$ of
$GL_2(\mathbf{A}_E)$ which are special or supercuspidal at all
places in $S(D)$.\end{thm}
\begin{proof}\cite{jaclan}.
\end{proof}
Let $\pi$ be an irreducible automorphic representation of
$G'(\mathbf{A}_\Q)$, and let $x$ be a place of $\Q$ which splits
in $K$. Then as in section \ref{carayol} we have an isomorphism
$\mathbf{Q}_x\isoto E_y$ as an $E$-algebra, and an identification
$G'(\Q_x)\cong B_y^\times\times\Q_x^\times$. Accordingly, we have
a decomposition $\pi_x\isoto\pi_y\otimes\psi_{\pi,y^c}$, where $c$
denotes complex conjugation. Let $BC(\pi_x)$ be the representation
$$\pi_y\otimes\pi_{y^c}\otimes(\psi_{\pi,y^c}\circ
c)\otimes(\psi_{\pi,y}\circ c)$$ of \begin{align*}G'(E_x)&\cong
B_x^\times\times E_x^\times\\&\cong B_y^\times\times
B_{y^c}^\times\times E_y^\times\times E_{y^c}^\times.\end{align*}
There is also, in \cite{ht01} a definition of $BC(\pi_x)$ for all
but finitely many places of $\Q$ which are inert in $E$. We will
not need this definition, except to note that by our choice of $E$
there is in fact a definition for \emph{all} inert places, and
that the base change of an unramified representation is unramified
(this follows easily from the discussion on page 199 of
\cite{ht01}).

\begin{thm}\label{172}Suppose that $\pi$ is an irreducible automorphic
representation of $G'(\A)$ such that $\pi_\infty$ has weight 2. Then
there is a unique irreducible automorphic representation
$BC(\pi)=(\psi,\Pi)$ of
$\mathbf{A}^\times_K\times(D\otimes\A)^\times$ such that
\begin{enumerate}
  \item $\psi=\psi_\pi|_{\mathbf{A}_K^\times}^c$.
  \item If $x$ is a place of $\Q$  then
  $BC(\pi)_x=BC(\pi_x)$.
  \item $\Pi_\infty$ has weight 2.
  \item $\psi_\Pi|_{\mathbf{A}_K^\times}=\psi^c/\psi$.
  \item $\Pi^*\cong\Pi$, where $\Pi^*(g)=\Pi(g^{-*})$.
\end{enumerate}
\end{thm}
\begin{proof}
  Immediate from Theorem VI.2.1 of \cite{ht01}.
\end{proof}
\subsection{Galois representations}\label{galrep} We recall a
version of Matsushima's formula (see page 8 of \cite{ht01} and
page 420 of \cite{car862}). For $i=0,1,2$ we have
$$\operatorname{H}^i_{\acute{e}t}(M'_{K'}\times\overline{E},\overline{\Q}_p)=\oplus_\pi\pi^{K'}\otimes
R^i(\pi),$$ where the sum is over irreducible admissible
representations $\pi$ of $G'(\A)$ with $\pi_\infty$ of weight 2, and
$R^i(\pi)$ is a certain finite dimensional continuous representation
of $\Gal(\overline{E}/E)$.

Define a virtual representation
$$[R(\pi)]=\sum_{i=0}^2(-1)^{i+1}[R^i(\pi)].$$ It can be shown that
$$[R(\pi)]\neq 0$$ (see the remarks after Theorem 1 in
\cite{ko92}). In the cases of interest to us,
$R^0(\pi)=R^2(\pi)=0$; in fact we have

\begin{thm}
  Let $\pi$ be an irreducible admissible representation of
  $G'(\A)$ over $\overline{\Q}_p$ with $\pi_\infty$ of weight 2, and
  suppose that $BC(\pi)=(\psi,\Pi)$ with $JL(\Pi)$
  cuspidal. Then $R^0(\pi)=R^2(\pi)=0$.
\end{thm}
\begin{proof}
  This is a special case of Corollary VI.2.7 of \cite{ht01}.
\end{proof}
Let $m_\pi$ denote the multiplicity of $\pi$ in the space of
automorphic forms on $G'(\Q)\backslash G'(\A)$ transforming by
$\psi_\pi$ under the centre of $G'(\A)$. Then
\begin{lemma}$\dim[R(\pi)]=2m_\pi=2$.\label{multone}
\end{lemma}
\begin{proof}
  The first equality follows at once from Theorem 1.3.1 of \cite{har00} and the
  remarks after Theorem 1 in \cite{ko92}. The statement that
  $m_\pi=1$ may be found on page 132 of \cite{cl99}.
\end{proof}

Thus under the assumption that $BC(\pi)=(\psi,\Pi)$ with $JL(\Pi)$
cuspidal, we have attached a continuous 2-dimensional $p$-adic
representation $R^1(\pi)$ to $\pi$ (note that we have implicitly
used the theory of the conductor for $U(1,1)$ developed in section
5 of \cite{lr04}).

By continuity $R^1(\pi)$ is defined over a finite extension
$K_\pi$ of $\Q_p$, and (from the compactness of
$\Gal(\overline{E}/E)$) there is a $\Gal(\overline{E}/E)$-stable
lattice $L$ in $R^1(\pi)$. Denote the semisimplification of the
representation $L/m_{K_\pi}L$ (where $m_{K_\pi}$ is the maximal
ideal of $\mathcal{O}_{K_\pi}$) by $\overline{\rho}_\pi$. It is
not \emph{a priori} obvious that $\overline{\rho}_\pi$ is
independent of the choice of $L$, but this follows from Theorem
\ref{frob} below; by the Cebotarev density theorem and the
Brauer-Nesbitt theorem (recalling that $p>2$)
$\overline{\rho}_\pi$ is determined by the knowledge of
$\tr(\overline{\rho}_\pi(\Frob_y))$ for all but finitely many $y$
lying above primes which split in $K$.

Before stating Theorem \ref{frob} it is convenient to recall the
connection between holomorphic differentials on $M'_{K'}$ and the
first \'{e}tale cohomology of $M'_{K'}$. We have the
Hodge-theoretic decomposition
\begin{align*}\operatorname{H}^0(M'_{K'}\otimes\C,\Omega^1_{M'_{K'}})\oplus\operatorname{H}^0(M'_{K'}\otimes\C,\overline{\Omega}^1_{M'_{K'}})&\isoto\operatorname{H}^1(M'_{K'}\otimes\C,\C)\\
&\isoto\operatorname{H}^1_{\acute{e}t}(M'_{K'}\times\overline{E},\overline{\Q}_p),\end{align*}from
which a standard argument shows that an automorphic form $\pi$ as
above corresponds to a unique differential $\omega_\pi$ on
$M'_{K'}$, which is an eigenvector for the Hecke operators
$T_\mathfrak{l}$. Specifically, this correspondence is determined
by the requirement that at all places $\mathfrak{l}$ of $E$ at
which $\pi_\mathfrak{l}$ is unramified principal series and at
which we have defined $T_\mathfrak{l}$, the eigenvalue of
$T_\mathfrak{l}$ acting on $\omega_\pi$ is $\mathbf{N}\l^{-1/2}$
multiplied by the sum of the Satake parameters of
$\pi_\mathfrak{l}$ (see for example \S 1.5 of \cite{ht02}). Then
we have
\begin{thm}\label{frob}Suppose that $K'=U_1(\mathfrak{n})$,
$y\nmid \mathfrak{n}p$ is a prime of $E$ lying over a place of
$\Q$ which splits in $K$, and $\pi_y$ is unramified principal
series. If $T_y\omega_\pi=a_y\omega_\pi$ and
$S_y\omega_\pi=b_y\omega_\pi$ then
$$\tr(\overline{\rho}_\pi(\Frob_y))=\overline{a_y/\psi_\pi(\Frob_y)}$$and
$$\det(\overline{\rho}_\pi(\Frob_y))=\overline{b_y/\psi_\pi(\Frob_y)},$$where
$\psi_\pi$ is the central character of $\pi$.
\end{thm}
\begin{proof}This is immediate from Corollary VII.1.10 of
\cite{ht01}.
\end{proof}
\begin{corollary}\label{154}If $\omega_f$ is a differential on
${\mathbf{M}'}^\text{can}_{\operatorname{\operatorname{bal}}.U_1(\p),H'}$ which is an eigenvector
for the operators $T_\mathfrak{l}$ with eigenvalues
$a_\mathfrak{l}$, then there is a continuous representation
$$\overline{\rho}_f:\Gal(\overline{E}/E)\to\GL_2(\overline{\mathbf{F}}_p)$$such
that for all but finitely many $\mathfrak{l}$ lying over primes of
$\Q$ which split in $K$ we have
$$\tr(\overline{\rho}_f(\Frob_\mathfrak{l}))=a_\mathfrak{l}.$$
\end{corollary}
\begin{proof}
  By the usual Deligne-Serre minimal prime lemma (Lemma 6.11 of
  \cite{delser}) there is a
  differential $\omega_\pi$ whose Hecke eigenvalues lift those of $\omega_f$; then the required
  representation is obtained by twisting the definition of
  $\overline{\rho}_\pi$ by $\overline{\psi}_\pi$.
\end{proof}
\section{Construction of Companion Forms}\label{construct}
\subsection{Hecke operators}

We define Hecke operators on $M'_{K'}$ in the
familiar way (namely as double cosets). We then give a modular
interpretation of these operators, which allows us to extend their
definitions to the integral models. Indeed, let $\mathfrak{l}$ be
a prime of $F$ lying over the prime $l\neq p$ of $\Q$, and suppose
that $B$ is split at $\mathfrak{l}$. Suppose also that $l$ splits
in $\Q(\sqrt{\lambda})$. Then if
$\mathfrak{l}_1=\mathfrak{l},\mathfrak{l}_2,\dots,\mathfrak{l}_k$
are the primes of $F$ above $l$, we have
$G'(\Q_l)\isoto\Q_l^*\times\GL_2(F_\mathfrak{l})\times\GL_2(F_{\mathfrak{l}_2})\times\dots\times\GL_2(F_{\mathfrak{l}_k})$.
Suppose now that $K'$ is of the form $\prod_qK'_q$, and that
$K'_l=\Z_l\times\prod_{i=1}^k \GL_2(\bigO_{F,\mathfrak{l}_i})$.
Let $\varpi_\l$ be an element of $\A^\infty$ which is a
uniformiser at $\mathfrak{l}$ and $1$ everywhere else; then we
define
$$T_\mathfrak{l}=\left[K'\begin{pmatrix}\varpi_\l&0\\0&1\end{pmatrix}K'\right].$$ In order to compute the action of
$T_\l$ we express it as the ratio of two other Hecke operators
$X_\l$ and $Y_\l$; define $X_\l=[K'x_\l K']$ and $Y_\l=[K'y_\l
K']$ where with obvious notation we have
$$x_\l=\left(l^{-1},\Bigl(\begin{smallmatrix}1&0\\0&\varpi_\l^{-1}\end{smallmatrix}\Bigr),\Bigl(\begin{smallmatrix}\varpi_{\l_2}^{-1}&0\\0&\varpi_{\l_2}^{-1}\end{smallmatrix}\Bigr),\dots,\Bigl(\begin{smallmatrix}\varpi_{\l_k}^{-1}&0\\0&\varpi_{\l_k}^{-1}\end{smallmatrix}\Bigr)\right)$$
and
$$y_\l=\left(l^{-1},\Bigl(\begin{smallmatrix}\varpi_\l^{-1}&0\\0&\varpi_\l^{-1}\end{smallmatrix}\Bigr),\Bigl(\begin{smallmatrix}\varpi_{\l_2}^{-1}&0\\0&\varpi_{\l_2}^{-1}\end{smallmatrix}\Bigr),\dots,\Bigl(\begin{smallmatrix}\varpi_{\l_k}^{-1}&0\\0&\varpi_{\l_k}^{-1}\end{smallmatrix}\Bigr)\right)$$
and we have $T_\l=X_\l Y_\l^{-1}$.

In the usual way we have  induced correspondences $X_\mathfrak{l}$
and $Y_\l$ on $\mathbf{M}_{K'}$. In section 7.5 of \cite{car861}
Carayol computes the action of $G'(\A^\infty)$ on the inverse
limit $\mathbf{M}'=\lim\limits_{\longleftarrow\atop K'}\mathbf{M}'_{K'}$, which
allows us to compute $X_\mathfrak{l}$ and $Y_\l$ on
${\mathbf{M}'}^\text{can}_{\operatorname{\operatorname{bal}}.U_1(\p),H'}$. In the usual way we
have an induced action on
$\Omega^1_{{\mathbf{M}'}^\text{can}_{\operatorname{\operatorname{bal}}.U_1(\p),H'}/\bigO_\p}$,
and Carayol's results allow us to explicate this as
$$\omega|X_\mathfrak{l}(A,\iota,\lambda,\eta_\p,\overline{\eta}_\p,P,Q)\mapsto\sum_{\phi}\omega\left((\phi A,\phi \iota,\phi \lambda,\phi\eta_\p,\phi\overline{\eta}_\p,\phi P,\phi
Q)\right),$$where the sum is over a certain set of
$\mathbf{N}\l+1$ isogenies of degree $l^d$ (we do not need to know
the precise details of these isogenies). Similarly, the action of
$Y_\l$ is via a single isogeny of degree $l^d$.

\subsection{Igusa curves}\label{hecke}We now extend the Hecke
action to an action on meromorphic functions on Igusa curves, and
relate this action to the action on differentials on Shimura
curves defined above. The required notation was introduced in
Chapter 2.

Let $\pi:\mathbf{A}'\to\mathbf{M}'_{0,H'}$ be the universal
abelian variety. The $\bigO_{\mathbf{M}'_{0,H'}}$-module
$\pi_*\Omega^1_{\mathbf{A}'/\mathbf{M}'_{0,H'}}$ is an
$\bigO_D\otimes\Z_p$-module, and we put
$$\omega_{\mathbf{A}'/\mathbf{M}'_{0,H'}}=\left(\pi_*\Omega^1_{\mathbf{A}'/\mathbf{M}'_{0,H'}}\right)^{2,1}_1.$$
Since $\Lie^{2,1}_1(\mathbf{A}'_{0,H'})$ is locally free of rank
one we see that $\omega_{\mathbf{A}'/\mathbf{M}'_{0,H'}}$ is
locally free of rank one. Similarly (or by pullback from
$\omega_{\mathbf{A}'/\mathbf{M}'_{0,H'}}$ using the universal
property of $\pi:\mathbf{A}'\to\mathbf{M}'_{0,H'}$) we have a
locally free sheaf of rank one
$\omega_{(\mathbf{A}')^\vee/\mathbf{M}'_{0,H'}}$. Then we have
\begin{thm}The Kodaira-Spencer map gives a canonical isomorphism
$$\omega_{\mathbf{A}'/\mathbf{M}'_{0,H'}}\otimes\omega_{(\mathbf{A}')^\vee/\mathbf{M}'_{0,H'}}\isoto\Omega^1_{\mathbf{M}'_{0,H'}/\bigO_\p}.$$
\end{thm}
\begin{proof}
  This is Proposition 4.1 of \cite{kas04}.
\end{proof}
We have a morphism
$\pi_{Ig}:\mathbf{M}'_{Ig,H'}\to\mathbf{M}'_{0,H'}\otimes\kappa$ of
degree $(p-1)$ which is \'{e}tale over the ordinary locus and
totally ramified over the supersingular locus (the former follows
from a computation of the local rings, and the latter from the
uniqueness of Igusa structures at supersingular points). Let $ss$
denote the divisor of supersingular points on $\mathbf{M}'_{Ig,H'}$,
and let $s$ be the degree of this divisor.

\begin{prop}$2s=(p-1)\deg\Omega^1_{\mathbf{M}'_{0,H}/\kappa}$.
\end{prop}
\begin{proof}By flatness we have
$\deg\Omega^1_{\mathbf{M}'_{0,H}/\kappa}=2(g-1)$, where $g$ is the
genus of any geometric fibre of $\mathbf{M}'_{0,H}/\bigO_\p$, so it
suffices to prove that $s=(p-1)(g-1)$. We do this by computing the
genus $g_0$ of $\mathbf{M'}_{U_0(\p),H'}$ in two ways. In
characteristic zero the Riemann-Hurwitz formula gives
$g_0=(p+1)(g-1)+1$, whereas on the special fibre we have
$g_0=2g+s-1$, from the description of the special fibre of
$\mathbf{M}'_{U_0(\p),H'}$ in section 10 of \cite{jar99} (or rather
from the obvious modification of the argument of \cite{jar99} to our
case). Equating these expressions gives the result.
\end{proof}
Put
$\omega^+:=\pi_{Ig}^*(\omega_{\mathbf{A}'/\mathbf{M}'_{0,H'}}\otimes\kappa)$
and
$\omega^-:=\pi_{Ig}^*(\omega_{(\mathbf{A}')^\vee/\mathbf{M}'_{0,H'}}\otimes\kappa)$.
\begin{prop}$\omega^+$ and $\omega^-$ have degree $s$, and there
is a natural isomorphism
$\Omega^1_{\mathbf{M}'_{Ig,H'}}\cong\omega^+\otimes\omega^-((p-2)(ss))$.
\end{prop}
\begin{proof}The degree of the polarisation associated to
$\mathbf{A}'$ is prime to $p$, so it is \'{e}tale, and thus
induces a (non-canonical) isomorphism $\omega^+\isoto\omega^-$, so
$\deg\omega^+=\deg\omega^-$. Then
\begin{align*}2s&=(p-1)\deg\Omega^1_{\mathbf{M}'_{0,H'}/\kappa}\\&=2(p-1)\deg(\omega_{\mathbf{A}'/\mathbf{M}'_{0,H'}}\otimes\kappa)\end{align*}so
that
\begin{align*}\deg\omega_{\mathbf{M}'_{Ig,H'}}&=\deg(\pi_{Ig})\deg(\omega_{\mathbf{A}'/\mathbf{M}'_{0,H'}}\otimes\kappa)\\ &=(p-1)\deg(\omega_{\mathbf{A}'/\mathbf{M}'_{0,H'}}\otimes\kappa)\\& =s.\end{align*}Then by the Riemann-Hurwitz formula we have
\begin{align*}\Omega^1_{\mathbf{M}'_{Ig,H'}}&\cong\pi_{Ig}^*\Omega^1_{\mathbf{M}'_{0,H'}}((p-2)ss)\\&\cong\pi_{Ig}^*((\omega_{\mathbf{A}'/\mathbf{M}'_{0,H'}}\otimes\kappa)\otimes(\omega_{(\mathbf{A}')^\vee/\mathbf{M}'_{0,H'}}\otimes\kappa))((p-2)ss)\\&\cong\omega^+\otimes\omega^-((p-2)ss).\end{align*}
\end{proof}
We now define sections $a^+$ of $\omega^+$ and $a^-$ of
$\omega^-$; $a^+\otimes a^-$ will play the role of a ``Hasse
invariant''. Given the data of an abelian variety with
polarisation and level structure $(\pi:A\to
S,\iota,\lambda,\eta_\p,\overline{\eta}_\p)$, where $S$ is an
$\F_p$-scheme, together with an Igusa structure
$P\in\ker(V|_{A^{\sigma}})$, we need to construct an element of
$\operatorname{H}^0(A,\Omega^1_{A/S})$.

By Cartier duality, $P$ gives a map
$\phi_P:\ker(F|_A)\to\mathbf{G}_m$. The standard invariant
differential $dX/X$ on $\mathbf{G}_m$ pulls back to an invariant
differential $\phi_P^*(dX/X)$ on $\ker(F|_A)$. Because $S$ is an
$\F_p$-scheme, the restriction map $$(\text{invariant 1-forms on }
A/S)\To(\text{invariant 1-forms on }\ker(F|_A))$$is an
isomorphism, so there is a unique invariant differential on $A$
whose restriction to $\ker(F|_A)$ is $\phi_P^*(dX/X)$.

Thus we have an element of
$\operatorname{H}^0(A,\Omega^1_{A/S})=\operatorname{H}^0(S,\pi_*\Omega^1_{A/S})$,
so an element of
$\operatorname{H}^0(S,(\pi_*\Omega^1_{A/S})^{2,1}_1)$. Applying
this with $S=\mathbf{M}'_{Ig,H'}$ and $A=\mathbf{A'}/S$ (the
universal abelian variety defined above) gives the required
section $a^+$ of $\omega^+$. The section $a^-$ of $\omega^-$ is
defined in the same way.
\begin{prop}$a^+$ has simple zeroes at the supersingular
points, and is nonzero on the ordinary locus.
\end{prop}
\begin{proof}If $A/\overline{\kappa}$ is supersingular, then the
only Igusa structure $P\in\ker(F|_A)$ is $P=0$, so $\phi_P=0$ and
we see that $a^+$ vanishes at every supersingular point. On the
other hand if $A/\overline{\kappa}$ is ordinary then we obviously
obtain a nonzero element of
$\operatorname{H}^0(S,(\pi_*\Omega^1_{A/S})^{2,1}_1)$, so $a^+$
cannot be identically zero. Since the degree of $\omega^+$ is $s$,
the number of supersingular points, $a^+$ can only have simple
zeroes at the supersingular points, and is nonzero
elsewhere.\end{proof} We define a Hecke action on
$\omega^+\otimes\omega^-$ by demanding that it commutes with the
Kodaira-Spencer isomorphism and the action on differentials
defined above.
\begin{prop}If $\nu\otimes\nu'$ is a local section of
$\omega^+\otimes\omega^-$, then
$$((\nu\otimes\nu')X_\mathfrak{l})(A,\iota,\lambda,\eta_\p,\overline{\eta}_\p,P)=\frac{1}{l^d}\sum_\phi\phi^*\nu|_{\phi(A,\iota,\lambda,\eta_\p,\overline{\eta}_\p,P)}\otimes\phi^*\nu'|_{\phi(A,\iota,\lambda,\eta_\p,\overline{\eta}_\p,P)}$$
where the sum is over the same isogenies used in the explication
of $X_\mathfrak{l}$ above, and a similar result holds for $Y_\l$.
\end{prop}
\begin{proof}We use the following description of the
Kodaira-Spencer map: for $A/S$ an abelian scheme, there is a short
exact sequence
$$0\to\Omega^1_{A/S}\to\operatorname{H}^1_{dR}(A/S)\to(\Omega^1_{A^\vee/S})^\vee\to
0,$$and a cup-product pairing
$$\langle,\rangle_{dR}:\operatorname{H}^1_{dR}(A/S)\times\operatorname{H}^1_{dR}(A/S)\to\mathcal{O}_S.$$If
$S$ is a $\Sigma$-scheme for some scheme $\Sigma$, we also have
the Gauss-Manin connection
$$\nabla:\operatorname{H}^1_{dR}(A/S)\to\operatorname{H}^1_{dR}(A/S)\otimes\Omega^1_{S/\Sigma}$$(see section III.9 of \cite{cf90}). This gives a map
\begin{align*}\omega^+\otimes\omega^-&\to\Omega^1_{S/\Sigma}\\(v\otimes
w)&\mapsto\langle v,\nabla w\rangle.\end{align*}Thus under the
Kodaira-Spencer isomorphism the image of
$$\phi^*\nu|_{\phi(A,\iota,\lambda,\eta_\p,\overline{\eta}_\p,P)}\otimes\phi^*\nu'|_{\phi(A,\iota,\lambda,\eta_\p,\overline{\eta}_\p,P)}$$
is the differential
$$\frac{1}{l^d}\sum_\phi\langle\phi^*\nu|_{\phi(A,\iota,\lambda,\eta_\p,\overline{\eta}_\p,P)},\nabla\phi^*\nu'|_{\phi(A,\iota,\lambda,\eta_\p,\overline{\eta}_\p,P)}\rangle.$$
Now, by the naturality of the Gauss-Manin connection we see that
if $\phi:A\to A'$ is an isogeny and $v'$ is an invariant
differential on $A'$, then $\nabla_A\phi^*v'=\phi^*\nabla_{A'}v'$
in $\operatorname{H}^1_{dR}$. Also, if $v$ is an invariant
differential on $A$, the fact that the adjoint of $A$ with respect
to the cup product pairing on de Rham cohomology is $\phi^\vee$,
satisfying $\phi^\vee\circ\phi=\deg\phi$, gives
$$\langle\phi^*v,\phi^*\nabla v'\rangle^A_{dR}=\deg\phi\langle
v,\nabla v'\rangle^{A'}_{dR}.$$ Every isogeny we sum over has
degree $l^d$, so we obtain (with obvious suppression of additional
structures)
\begin{align*}\frac{1}{l^d}\sum_\phi\langle\phi^*\nu|_{\phi(A)},\nabla\phi^*\nu'|_{\phi(A)}\rangle&=\sum_\phi\langle\nu|_{\phi(A)},\nabla\nu'|_{\phi(A)}\rangle\\
&=X_\mathfrak{l}(\langle\nu|_{(A)},\nabla\nu'|_{(A)}\rangle).\end{align*}
\end{proof}
We now examine two ways of moving between meromorphic functions
and meromorphic differentials on the Igusa curve. The first is the
map $d:f\mapsto df$, and the second is multiplication by
$a^+\otimes a^-$ followed by application of the Kodaira-Spencer
isomorphism.

We define a Hecke action on meromorphic functions on
$\mathbf{M}'_{Ig,H'}$ by
$$(X_\mathfrak{l}f)(A,\iota,\lambda,\eta_\p,\overline{\eta}_\p,P)=\frac{1}{l^d}\sum_\phi\phi^*f(\phi(A,\iota,\lambda,\eta_\p,\overline{\eta}_\p,P)),$$and
similarly for $Y_\l$.
\begin{thm}\label{126}$X_\mathfrak{l}(df)=l^d d(X_\mathfrak{l}f)$, and $Y_\mathfrak{l}(df)=l^d d(Y_\mathfrak{l}f)$, so that $T_\mathfrak{l}(df)= d(T_\mathfrak{l}f)$.
\end{thm}
\begin{proof}We have\begin{align*}l^d d(X_\mathfrak{l}f)&=l^d\frac{1}{l^d}d\left(\sum_\phi\phi^*f(\phi(A,\iota,\lambda,\eta_\p,\overline{\eta}_\p,P))\right)\\
&=\sum_\phi\phi^*(df)(\phi(A,\iota,\lambda,\eta_\p,\overline{\eta}_\p,P))\\
&=X_\mathfrak{l}(df).\end{align*}The result for $Y_\l$ follows
similarly.
\end{proof}
\begin{thm}$X_\mathfrak{l}((a^+\otimes a^-)f)=(a^+\otimes
a^-)X_\mathfrak{l}f$ and $Y_\mathfrak{l}((a^+\otimes
a^-)f)=(a^+\otimes a^-)Y_\mathfrak{l}f$, so that
$T_\mathfrak{l}((a^+\otimes a^-)f)=(a^+\otimes
a^-)T_\mathfrak{l}f$.
\end{thm}
\begin{proof}We compute\begin{align*}X_\mathfrak{l}((a^+\otimes a^-)f)&=\frac{1}{l^d}\sum_\phi\phi^*((a^+\otimes
a^-)f)(\phi(A,\iota,\lambda,\eta_\p,\overline{\eta}_\p,P))\\
&=\frac{1}{l^d}\sum_\phi\phi^*(a^+)|_{\phi
(A,P)}\otimes\phi^*(a^-)|_{\phi (A,P)}\phi^*f|_{\phi
(A,P)}\end{align*} (with obvious suppression of additional
structures other than the Igusa structure), and similarly for
$Y_\l$. It is thus sufficient to prove that for all isogenies
$\phi:(A,P)\to (A',Q)$ we have
$\phi^*(a^+|_{(A',Q)})=a^+|_{(A,P)}$ and
$\phi^*(a^-|_{(A',Q)})=a^-|_{(A,P)}$. We prove the former, the
proof of the latter being formally identical. We have maps
$\psi_P:\ker(F|_A)\to\mathbf{G}_m$,
$\psi_Q:\ker(F|_{A'})\to\mathbf{G}_m$, and
$a^+|_{\ker(F|_A)}=\psi_P^*(dX/X)$,
$a^+|_{\ker(F|_{A'})}=\psi_Q^*(dX/X)$. Since $a^+|_A$ is
determined by $a^+|_{\ker(F|_A)}$, it suffices to check that
$\psi_Q\circ\phi=\psi_P$, which is obvious.
\end{proof} We now examine an analogue of Atkin-Lehner's $U$
operator. Define, in the same fashion as $T_\mathfrak{l}$, an
operator
$$U_\p=\left[K'\begin{pmatrix}\varpi_\p&0\\0&1\end{pmatrix}K'\right]$$
on $\mathbf{M}'_{K'}$, where $\varpi_\p$ is a finite adele which is
a uniformiser of $\bigO_\p$ at $\p$ and $1$ everywhere else. This
acts as a correspondence on
${\mathbf{M}'}^\text{can}_{\operatorname{\operatorname{bal}}.U_1(\p),H'}$, and by the
functoriality of N\'{e}ron models we get an induced endomorphism of
$\Pic^0({\mathbf{M}'}^\text{can}_{\operatorname{\operatorname{bal}}.U_1(\p),H'}\otimes\overline{\kappa})$,
and in particular of $\Pic^0(\mathbf{M}'_{Ig,H'})$. This action is
not, however, given by a correspondence. Let $\Frob$ be the
Frobenius correspondence on $\mathbf{M}'_{Ig,H'}$, with dual $\Ver$.
Then we have \begin{thm}\label{Up}$U_\p=\Ver$ on
$\mathbf{M}'_{Ig,H'}$.
\end{thm}

\begin{proof}The proof is very similar to that of the
corresponding result in the classical case, for which see theorem
5.3 of \cite{wil80} or page 255 of \cite{mw84}.

One expresses $U_\p$ as a ratio of two Hecke operators which act
via isogenies, and computes on the ordinary locus.
\end{proof}
\subsection{Computations of cohomology classes}\label{section11}Suppose now that $k\geq
3$. Recall that we have associated an automorphic form $\pi_1$ on
$G'(\mathbf{A}_\Q)$ to our Hilbert modular form $\pi$, together
with a Galois representation $\overline{\rho}_f$. We write
$\epsilon'$ for the prime-to-$p$ part of the central character of
$\pi_1$.

We now associate a regular differential (cf. section 8 of
\cite{gro90}) $\omega_f$ on $\mathbf{M}'_{K'}$ to $\pi_1$, as in
section \ref{galrep}. Write $a_\p$ for the eigenvalue of $U_\p$ on
$\omega_f$. Let $\omega_A$ denote the differential on
$\mathbf{M}'_{Ig,H'}$ corresponding to $a^+\otimes a^-$ under the
Kodaira-Spencer isomorphism. We have an obvious action of
$(\bigO_\p/\p)^{*}$ on $\mathbf{M}'_{Ig,H'}$, which we denote by
$\langle\cdot\rangle$, so that
$\langle\alpha\rangle\omega_f=\alpha^{k-2}\omega_f=\alpha^{-k'}\omega_f$,
where $k':=p+1-k$. From the definitions of $a^+$ and $a^-$ we see
that $\langle\alpha\rangle a^+=\alpha^{-1}a^+$ and
$\langle\alpha\rangle a^-=\alpha^{-1}a^-$, so
$\langle\alpha\rangle\omega_A=\alpha^{-2}\omega_A$.

Let $v$ be a local parameter at a fixed supersingular point $y$ on
$\mathbf{M}'_{Ig,H'}$ such that $\langle\alpha\rangle
v=\alpha^{-1}v$, so we may expand
\begin{align*}\omega_A&=\left(\sum_{n=1}^\infty
c_n(v)v^{2+n(p-1)}\right)\frac{dv}{v}\\
\omega_f&=\left(\sum_{n=0}^\infty
a_n(v)v^{k'+n(p-1)}\right)\frac{dv}{v}.\end{align*}Define an
operator $M$ on meromorphic differentials on $\mathbf{M}'_{Ig,H'}$
by
$$M\omega:=d\left(\frac{\omega}{\omega_A}\right).$$Then we have
\begin{thm}\label{121}$\frac{M^{k'-1}\omega_f}{\omega_A}=\frac{(k'-1)!}{c_1(v)^{k'}}\left(a_0(v)v^{-pk'}+a_{k'}(v)v^{-k'}+\dots\right)$
\end{thm}
\begin{proof}Firstly, we note that $c_2(v)=0$. This follows
because, by Theorem \ref{Up}, we have
\begin{align*}\left(\sum_{n=1}^\infty
c_n(v)v^{2+n(p-1)}\right)\frac{dv}{v}&=\omega_A\\&=U_\p\omega_A\\&=\left(\sum_{n=0}^\infty
c_{pn+2}(v^\sigma)v^{2+n(p-1)}\right)\frac{dv}{v}.\end{align*}It
is easy to see that the leading term of the right hand side of the
proposed equality is correct, and this reduces us to computing
$M^{k'}\omega_f$. If we use the relation
$\frac{d}{dv}(v^pg)=v^p\frac{dg}{dv}$, we can turn the problem
into one in characteristic zero; if we define formal power series
\begin{align*}\omega'_A&=\left(\sum_{n=1}^\infty
C_n(v)v^{2-n}\right)\frac{dv}{v}\\
\omega'_f&=\left(\sum_{n=0}^\infty
A_n(v)v^{k'-n}\right)\frac{dv}{v},\end{align*}with $C_2(v)=0$ and
an operation $M'$ on differentials by
$$M'\omega:=d\left(\frac{\omega}{\omega'_A}\right),$$then it is enough to check that $$M'^{k'}\omega'_f=\left(-\frac{k'!}{C_1(v)^{k'}}A_{k'}(v)v^{-k'}+\dots\right)\frac{dv}{v}.$$ We can easily
convert our problem into one about power series, rather than
differentials; define
\begin{align*}D&=\left(\sum_{n=1}^\infty
C_n(v)v^{1-n}\right)^{-1}\\&=\frac{1}{vC_1(v)}\left(\sum_{n=0}^\infty
D_n(v)v^{-n}\right),\end{align*}say, with $D_1(v)=0$, and
$$f=\sum_{n=0}^\infty
A_n(v)v^{k'-1-n}.$$ Then if we define an operator $A$ on power
series by $Ag=vD\frac{d}{dv}g$, it is not hard to see that it
suffices to prove that the power series $D^{-1}A^{k'}Df$ has no
terms in $v^{-1},\dots,v^{1-k'}$. Then it suffices to prove that
$A^{k'}v^j$ has no terms in $v^{-1},\dots,v^{-k'}$ if $j\leq k'$.
We claim that we can write $v=cw+O(1)$ so that $Af=df/dw$; this is
equivalent to $dw/dv=1/(vD)$, and the integrability of $1/(vD)$
follows from its lack of a term in $v^{-1}$. Then as $v^j=O(w^j)$,
we see that the $k'$-th derivative of $v^j$ with respect to $w$ is
of the form $\mbox{constant}+O(w^{-k'+1})$ (where the constant is
zero unless $j=k'$), as required.
\end{proof}
\begin{corollary}\label{poles}$M^{k'}\omega_f=\left(-\frac{k'!}{c_1(v)^{k'}}a_{k'}(v)v^{-k'}+\dots\right)\frac{dv}{v}$.
\end{corollary}
\begin{proof}This is immediate\end{proof} This result, showing that $M^{k'}\omega_f$ has a pole of low
degree, allows us to associate a certain cohomology class to
$M^{k'}\omega_f$. Let $ss$ denote the divisor of supersingular
points in $\mathbf{M}'_{Ig,H'}$, and let
$U=\mathbf{M}'_{Ig,H'}\backslash ss$, and let $\eta$ denote the
generic point of $\mathbf{M}'_{Ig,H'}$.
\begin{prop}
The complex of groups
\begin{multline*}\mathcal{O}_{\mathbf{M}'_{Ig,H'}}(U)\to\Omega^1_{\mathbf{M}'_{Ig,H'}/F_\p^p}\oplus\mathcal{O}_{\mathbf{M}'_{Ig,H'},\eta}/\mathcal{O}_{\mathbf{M}'_{Ig,H'},ss}\\ \to
\Omega^1_{\mathbf{M}'_{Ig,H'},\eta}/\Omega^1_{\mathbf{M}'_{Ig,H'}/F_\p^p}(\log
ss)_{ss},\end{multline*}where the first arrow takes $h$ to $(dh,h)$,
and the second arrow takes a pair $(\omega,g)$ to $\omega-dg$,
computes the de Rham cohomology of $\mathbf{M}'_{Ig,H'}$ with log
poles on $ss$.
\end{prop}
\begin{proof}\cite{cv92}, \S2.
\end{proof}
Now, $M^{k'}\omega_f$ has poles of order at most $k'+1=p+2-k\leq
p-1$ on $ss$, by Corollary \ref{poles}, so there is a section $h$ of
$\mathcal{O}_{\mathbf{M}'_{Ig,H'}}((p-1)ss)_{ss}$ such that
$M^{k'}\omega_f-dh$ has at worst simple poles on $ss$. Such an $h$
is well defined modulo $\mathcal{O}_{\mathbf{M}'_{Ig,H'},ss}$, so we
have a well-defined cohomology class
$[M^{k'}\omega_f]\in\mathbf{H}^1(\mathbf{M}'_{Ig,H'},\Omega^\bullet_{\mathbf{M}'_{Ig,H'}}(\log
ss))$. Furthermore, since $k\neq 2$ Corollary \ref{poles} shows that
$M^{k'}\omega_f$ has zero residues, so in fact
$[M^{k'}\omega_f]\in\mathbf{H}^1(\mathbf{M}'_{Ig,H'},\Omega^\bullet_{\mathbf{M}'_{Ig,H'}})$.

Coleman has shown that there is an isomorphism between
$\mathbf{H}^1(\mathbf{M}'_{Ig,H'},\Omega^\bullet_{\mathbf{M}'_{Ig,H'}})$
and the quotient of the space of meromorphic differentials on
$\mathbf{M}'_{Ig,H'}$ with no residues and poles of order at most
$p$ on $\mathbf{M}'_{Ig,H'}$ by the space of exact differentials
$dg$, where $g$ is a meromorphic function on $\mathbf{M}'_{Ig,H'}$
with poles of order at most $(p-1)$ (for a statement, see page 499
of \cite{gro90}; the proof is a straightforward exercise in \v{C}ech
cohomology). Thus if $[M^{k'}\omega_f]=0$ there is a meromorphic
function $h$ on $\mathbf{M}'_{Ig,H'}$ with poles only at the
supersingular points, of order at most $(p-1)$, satisfying
$dh=M^{k'}\omega_f$. Furthermore, replacing $h$ by an element of the
vector space spanned by $\langle h,U_\p h,\dots\rangle$, we see that
there is a meromorphic function $h$ with poles of order at most
$(p-2)$ on $\mathbf{M}'_{Ig,H'}$, satisfying $U_\p h=b_\p h$ for
some $b_\p$.

\begin{lemma}
$T_\mathfrak{l}h=a_\mathfrak{l}h$.
\end{lemma}
\begin{proof}We have $dh=M^{k'}\omega_f$, so
$M(h\omega_A)=M^{k'}\omega_f$. Thus by Theorem \ref{126} we have
have \begin{align*}
M(T_\l(h\omega_A))&=T_\l(M(h\omega_A))\\&=T_\l(M^{k'}\omega_f)\\&=M^{k'}(T_\l\omega_f)
\\&=a_\l M^{k'}\omega_f\\&=a_\l
M(h\omega_A)\end{align*}so $M(T_\l(h\omega_A)-a_\l h\omega_A)=0$,
and $T_\l(h\omega_A)-a_\l h\omega_A=g^p\omega_A$ for some $g$. But
the equation $T_\l h-a_\l h =g^p$ immediately implies that $g=0$
unless $k=p$, upon comparison of leading terms. If $k=p$, then we
use the fact that $T_\l h-a_\l h$ is an eigenform for $U_\p$, with
eigenvalue $b_\p$. Then $b_\p g^p=U_\p g^p=g$, which again gives
$g=0$.\end{proof}

From Corollary \ref{poles} we have $$\langle\alpha\rangle
h=\alpha^{-k'}h.$$ Then the holomorphic differential
$\omega_{f'}=\omega_A h$ satisfies
\begin{align*}T_\mathfrak{l}\omega_{f'}&=a_\mathfrak{l}\omega_{f'}\\
\langle\alpha\rangle\omega_{f'}&=\alpha^{-k'}\omega_{f'}\\
U_\p\omega_{f'}&=b_\p\omega_{f'}\end{align*}for some $b_\p$.

Thus we have:
\begin{thm}\label{313}If $[M^{k'}\omega_f]=0$, then there is an ordinary ``companion
form'' $\omega_{f'}$. The representation $\overline{\rho}_{f'}$
attached to $\omega_{f'}$ is isomorphic to
$\overline{\rho}_{f}$.\end{thm}
\begin{proof}We have everything except for the assertion that
$\overline{\rho}_{f'}\cong\overline{\rho}_f$, and the claim that
$b_\p\neq 0$. The first follows at once from the Cebotarev density
theorem and the fact that
$T_\mathfrak{l}\omega_{f'}=a_\mathfrak{l}\omega_{f'}$. To prove
that $b_\p\neq 0$, note simply that if $b_\p=0$ then
$\omega_{f'}=M^{k'-1}\omega_f$, which is nonsense (from a
comparison of leading terms).
\end{proof}
We must now show that if $\overline{\rho}_{f}$ is unramified at
$\p$, then $[M^{k'}\omega_f]=0$.
\begin{thm}\label{unitroot}$\sigma^{-1}\Frob[M^{k'}\omega_f]=\frac{a_\p}{\epsilon(\p)}[M^{k'}\omega_f].$\end{thm}
\begin{proof}The proof of this is similar to the proof of Theorem
5.1 of \cite{cv92} (note however that our proof works for all
$3\leq k \leq p$, and may be used to replace the appeal to rigid
analysis in \cite{cv92}).

By Corollary \ref{poles}, $[M^{k'}\omega_f]$ is represented by the
cocycle
$$\left(M^{k'}\omega_f,\frac{(k'-1)!}{c_1(v)^{k'}}a_{k'}(v)v^{-k'}\right).$$
Since $M^{k'}\omega_f=d(M^{k'-1}\omega_f/\omega_A)$, this class is
also represented by
$$\left(0,-\frac{(k'-1)!}{c_1(v)^{k'}}a_{0}(v)v^{-pk'}\right).$$

Then if $\phi$ is the Frobenius endomorphism of
$\mathbf{M}'_{Ig,H'}$, $\sigma^{-1}\Frob[M^{k'}\omega_f]$ is
represented by
$$\phi^*\left(M^{k'}\omega_f,\frac{(k'-1)!}{c_1(v)^{k'}}a_{k'}(v)v^{-k'}\right)=\left(0,\frac{(k'-1)!}{c_1(v^\sigma)^{k'}}a_{k'}(v^\sigma)v^{-pk'}\right).$$We
thus need to demonstrate the equality
$$a_{k'}(v^\sigma)c_1(v)^{k'}=-\frac{a_\p}{\epsilon(\p)}a_0(v)c_1(v^\sigma)^{k'}.$$The equality
$U_\p\omega_f=a_\p\omega_f$ yields $a_{k'}(v)=a_\p a_0(v^\sigma)$
and thus $a_{k'}(v^\sigma)=a_\p a_0(v^{\sigma^2})$. A
straightforward computation shows that on the supersingular locus
we have the equality
$\sigma^{-2}=\langle\p\rangle\cdot\langle-1\rangle_\p$, where
$\langle\p\rangle$ is the Hecke operator
$\bigl(\begin{smallmatrix}\p&0\\0&\p\end{smallmatrix}\bigr)$ and
$\langle -1\rangle_\p$ is given by
$$ \langle -1\rangle_\p:(A,\iota,\lambda,\eta_\p,\overline{\eta}_\p,P,Q)\mapsto(A,\iota,\lambda,\eta_\p,\overline{\eta}_\p,-P,-Q).$$ We thus have $(\omega_f)^{\sigma^{-2}}=-\epsilon(\p)\omega_f$, from
which we obtain
$a_0(v^{\sigma^2})(v^{\sigma^2}|\langle\p\rangle^{-1}\langle
-1\rangle_\p^{-1})^{k'}=-\frac{1}{\epsilon(\p)}a_0(v)v^{k'}$.
Combining these results, we
have\begin{align*}a_{k'}(v^\sigma)c_1(v)^{k'}&=a_\p
a_0(v^{\sigma^2})v^{k'}c_1(v)^{k'}\\&=-\frac{a_\p}{\epsilon(\p)}a_0(v)\left(\left(\frac{v}{v^{\sigma^2}|\langle\p\rangle^{-1}\langle
-1\rangle_\p^{-1}}\right) (y)c_1(v)\right)^{k'},\end{align*}and
the result follows from Lemma \ref{B2}.\end{proof}

\begin{lemma}\label{B2}$\left(\frac{v}{v^{\sigma^2}|\langle\p\rangle^{-1}\langle
-1\rangle_\p^{-1}}\right)(y)=c_1(v^\sigma)/c_1(v)$.\end{lemma}
\begin{proof}This may be
proved in exactly the same fashion as Lemma 5.3 of \cite{cv92},
except that rather than working with the canonical elliptic curve
over $\F_{p^2}$ one works with the canonical formal
$\bigO_\p$-module defined over $\F_{q^2}$, for which see
Proposition 1.7 of \cite{dri76}.
\end{proof}

\begin{defn}As in Corollary \ref{154}, let $\omega_F$ be a regular differential on ${\mathbf{M}'}^\text{can}_{\operatorname{\operatorname{bal}}.U_1(\p),H'}$
which is an eigenform for the Hecke operators $T_\l$, $U_\p$, and
$\langle\alpha\rangle$, with eigenvalues lifting those of
$\omega_f$. We also write $\omega_F$ for the corresponding
differential on
${\mathbf{M}'}^\text{can}_{\operatorname{\operatorname{bal}}.U_1(\p),H'}$.\end{defn}
\begin{defn}Define a Hecke operator $w_\p$ on
${\mathbf{M}'}^\text{can}_{\operatorname{\operatorname{bal}}.U_1(\p),H'}$ by
$\bigl(\begin{smallmatrix}0&1\\-\p&0\end{smallmatrix}\bigr)$. This
gives a Hecke operator $w_\p$ on $\mathbf{M}'_{\operatorname{\operatorname{bal}}.U^p_1(\p),H'}$
by descent.\end{defn}

\begin{thm}\label{Upw}Let $X'$ be the subscheme of $\mathbf{M}'_{\operatorname{bal}.U^p_1(\p),H'}$
obtained by removing the closed subscheme
$\mathbf{M}'^{(\sigma^{-1})}_{Ig,H'}$. Then $(\omega_F|U_\p
w_\p)|_{X'}$ is divisible by $\pi^{k'}$, so that there is a
cohomology class $[(\omega_F|U_\p w_\p)|_{X'}/\pi^{k'}]$ on
$\mathbf{M}'_{Ig,H'}$. Then $$[(\omega_F|U_\p
w_\p)|_{X'}/\pi^{k'}]=-u\frac{\epsilon(\p)}{k'!}[M^{k'}\omega_f],$$for
some unit $u$.
\end{thm}
\begin{proof}Recall that the completed local ring of
${\mathbf{M}'}^\text{can}_{\operatorname{bal}.U_1(\p),H'}$ at a supersingular point
is $\bigO_{F^0_\p}[[v,w]]/(vw-\pi)$, where
$\mathbf{M}'^{(\sigma^{-1})}_{Ig,H'}$ is given by $w=0$. Thus on
$X'$ the function $w$ is invertible, and we have $v=\pi/w$. An easy
check shows that on the supersingular locus we have
$w_\p=\sigma\cdot\langle\p\rangle$, whence we have

\begin{align*}(\omega_F|U_\p
w_\p)|_{X'}&=\left(\epsilon(\p)a_{k'}(w)(w^\sigma)^{k'}+\dots\right)\frac{dv}{v}\text{
mod
}p\\&=\left(u\epsilon(\p)a_{k'}(v)\frac{\pi^{k'}}{v^{k'}c_1(v)^{k'}}+\dots\right)\frac{dv}{v}\text{
mod }p,\end{align*}where we have used the identity
$c_1(v)^{-1}=uv(v^{\sigma}|w_\p\langle\p\rangle^{-1})$, which may
again be proved in exactly the same fashion as Lemma 5.3 of
\cite{cv92}, except that rather than working with the canonical
elliptic curve over $\F_{p^2}$ one works with the canonical formal
$\bigO_\p$-module defined over $\F_{q^2}$, for which see Proposition
1.7 of \cite{dri76}.

So $(\omega_F|U_\p w_\p)|_{X'}$ is divisible by $\pi^{k'}$, and
comparing this with Corollary \ref{poles} we see that the claim
follows.

\end{proof}

\subsection{Pairings}\label{12}We now recall from \cite{cv92} the
relationship between the Kodaira-Spencer and Serre-Tate pairings,
and a formula relating the Kodaira-Spencer pairing on a
semi-stable curve to the cup product of certain (log-) cohomology
classes.

Let $R$ be a complete local ring with residue field $\F$ of
characteristic $p$, and let $W(\F)$ denote the ring of Witt
vectors of $\F$. Let $G$ be a $p$-divisible group over $R$ with
dual $G^\vee$, and let $\Omega_G$, $\Omega_{G^\vee}$ denote the
invariant one-forms on $G/R$, $G^\vee/R$ respectively. Then (see
Corollary 4.8.iii of \cite{ill85}) there is a pairing, functorial
for morphisms of $p$-divisible groups over $R$,
$$\kappa:\Omega_G\otimes\Omega_{G^\vee}\to\Omega^1_{R/W(\F)}.$$Let
$\overline{G}$, $\overline{G}^\vee$ denote the special fibres of
$G$, $G^\vee$ respectively. We have the Serre-Tate pairing (see
\cite{kat81}) $$q:T_p\overline{G}\times T_p\overline{G}^\vee\to
1+m_R,$$where $T_p$ denotes the $p$-adic Tate module, and $m_R$ is
the maximal ideal of $R$. We can view an element $\alpha\in T_pG$
as a homomorphism from $G^\vee$ to $\mathbf{G}_m$, and we define
$\omega_\alpha=\alpha^*(dt/t)\in\Omega_{G^\vee}$. If
$\alpha^\vee\in T_p G^\vee$ we define
$\omega_{\alpha^\vee}\in\Omega_G$ in the same fashion. For $a\in
R^*$, let $d\log(a)=d_{R/W(\F)}a/a\in\Omega^1_{R/W(\F)}$; then we
have

\begin{thm}Suppose that $G$ is an ordinary $p$-divisible group
over $R$ (that is, suppose that the dual of the connected subgroup
of $\overline{G}$ is \'{e}tale). Then for any $\alpha\in T_pG$,
$\alpha^\vee\in T_pG^\vee$, we have $$d\log
q(\alpha,\alpha^\vee)=\kappa(\omega_{\alpha^\vee}\otimes\omega_\alpha).$$
\end{thm}
\begin{proof}This is Theorem 1.1 of \cite{cv92}.\end{proof}

Let $S=\F[t]/(t^{b+1})$ with $0\leq b<p$. Let $S^\times$ denote
the log-scheme associated to the pre-log structure $\mathbf{N}\to
S$, $1\mapsto t$. Let $M_S$ denote the corresponding monoid, with
an element $T$ mapping to $t$. Let $\F$ (slightly abusively)
denote $\F$ with the trivial log-structure; it follows easily that
$\Omega^1_{S^\times/\F}$ is a free $S$-module generated by $d\log
T$.

Let $s:X\to\Spec(R)$ be a semi-stable curve over $S$ (that is, $X$
is locally isomorphic to $xy=t$ in the \'{e}tale topology), and
suppose that there is a lifting $\tilde{X}$ of $X$ to a
semi-stable curve over $\tilde{S}$, where $\tilde{S}$ is a
discrete valuation ring with $\tilde{S}\mbox{ mod }p=S$ and the
generic fibre of $\tilde{X}$ is smooth over the generic point of
$\tilde{S}$ (this will, of course, hold in our applications to PEL
Shimura curves). We have natural log-structures on $\tilde{X}$ and
$\tilde{S}$ given by the subsheaf of the structure sheaf whose
sections become invertible upon removal of the special fibre. Let
$X^\times$ be the reduction of this log-scheme to $S$; $X^\times$
is smooth over $S^\times$.

We have an exact sequence of sheaves (see Proposition 3.12 of
\cite{kat89})
$$0\to s^*\Omega^1_{S^\times/\F}\to\Omega^1_{X^\times/\F}\to\Omega^1_{X^\times/S^\times}\to
0.$$Let
$Kod:\operatorname{H}^0(X,\Omega^1_{X^\times/S^\times})\to\operatorname{H}^1(X,s^*\Omega^1_{S^\times/\F})\cong\operatorname{H}^1(X,\bigO_X)\otimes
\Omega^1_{S^\times/\F}$ denote the boundary map in the
corresponding long exact sequence of cohomology. When $X$ is
smooth over $S$ there is a natural isomorphism
$\Omega^1_{X^\times/S^\times}\isoto\Omega^1_{X/S}$, and $Kod$ is
the composition of the usual Kodaira-Spencer map and the natural
map
$\operatorname{H}^1(X,\bigO_X)\otimes\Omega^1_{S/\F}\to\operatorname{H}^1(X,\bigO_X)\otimes\Omega^1_{S^\times/\F}$.

Suppose that the reduction of $X$ mod $t$ is $\overline{X}=C_1\cup
C_2$ with $C_1$, $C_2$ smooth irreducible curves. Let $D=C_1\cap
C_2$, $U_1=X-C_2$, $U_2=X-C_1$. Let $d$ denote the boundary map
for the complex $\Omega^\bullet_{X^\times/S^\times}$. Then we
have:

\begin{thm}Suppose $\omega$ is in the image of the natural map
from $\operatorname{H}^0(X,\Omega^1_{X/S})$ to
$\operatorname{H}^0(X,\Omega^1_{X^\times/S^\times})$ and that
$\omega|_{U_2}=t^b\eta$ for $\eta\in\Omega^1_{X/S}(U_2)$. Then
$$\overline{\eta}:=\eta|_{C_2}\in\operatorname{H}^0(C_2,\Omega^1_{C_2/\F}((b+1)D)).$$In
particular, we obtain a cohomology class
$$[\overline{\eta}]\in\mathbf{H}^1(C_2,\Omega^1_{C_2/\F}(\log D))$$
in the same fashion as in Section \ref{section11}. If $\nu$ is in
the image of $\operatorname{H}^0(X,\Omega^1_{X/S})$ in
$\operatorname{H}^0(X,\Omega^1_{X^\times/S^\times})$ and
$\nu|_{U_1}\in t^b\Omega^1_{X^\times/S^\times}(U_1)$, then
$$\nu\cdot Kod(\omega)=([\nu|_{C_2}],[\eta|_{C_2}])_{C_2}bt^bd\log
T,$$where $(,)_{C_2}$ is the usual cup product pairing.

\end{thm}
\begin{proof}This is Theorem 3.1 of \cite{cv92}.\end{proof}
\subsection{The local Galois Representation at $\p$}In this section we examine the local Galois representation
on the Jacobian of ${\mathbf{M}'}^\text{can}_{\operatorname{bal}.U_1(\p),H'}$. We
write $J:=\Jac({\mathbf{M}'}^\text{can}_{\operatorname{bal}.U_1(\p),H'})$,
$\mathbf{T}\subseteq\End_\Q(J)$ generated by the $U_\p$,
$T_\mathfrak{l}$, $S_\mathfrak{l}$ for $\mathfrak{l}$ unramified and
not in the level, and the operators $\langle \alpha\rangle$. Let
$\mathfrak{m}$ be a maximal ideal of $\mathbf{T}$ above the (minimal
prime) ideal corresponding to $\omega_f$. $\mathbf{T}$ is free of
finite rank over $\mathbf{Z}$, so the ring $\mathbf{T}_p=\varprojlim
\mathbf{T}/p^n\mathbf{T}=\mathbf{T}\otimes\mathbf{Z}_p$ is a
complete semilocal $\mathbf{Z}_p$-algebra of finite rank. The ring
$\mathbf{T}_\mathfrak{m}=\varprojlim\mathbf{T}/\mathfrak{m}^n\mathbf{T}$
is complete and local, and by the theory of complete semilocal rings
$\mathbf{T}_\mathfrak{m}$ is a direct factor of $\mathbf{T}_p$, so
we have an idempotent decomposition of the identity
\begin{align*}\mathbf{T}_p&=\mathbf{T}_\mathfrak{m}\times\mathbf{T}'_\mathfrak{m}\\1&=\epsilon_\mathfrak{m}+\epsilon_\mathfrak{m'}.\end{align*}
Let $h=\rk_{\mathbf{Z}_p}(\mathbf{T}_\mathfrak{m})$. Let $G$ be
the $p$-divisible group over $E$ defined by
$T_pG=\epsilon_\mathfrak{m}T_pJ$.

\begin{thm}\begin{enumerate}\item The $p$-divisible group $G$ has
height $2h$ and is isomorphic to $G^\vee$ over $E_\mathfrak{np}$,
where $E_\mathfrak{np}$ is the ray class field of conductor
$\mathfrak{np}$. It has good reduction over $F_\mathfrak{p}'$,
where $F_\mathfrak{p}'$ is the tame extension of $F_\p$ of degree
$q-1$ corresponding to our choice of uniformiser $\p$.
\item Let $\overline{G}$ be the reduction of $G$ over
$\bigO_{F_\mathfrak{p}'}/\mathfrak{m}_{\bigO_{F_\mathfrak{p}'}}\bigO_{F_\mathfrak{p}'}=\mathbf{F}_p$,
and $D(\overline{G})$ be its Dieudonn\'{e} module. Then
$\overline{G}=\overline{G}^m\times\overline{G}^e$, where
$\overline{G}^m$ is multiplicative and $\overline{G}^e$ is
\'{e}tale, and $\overline{G}^m$ and $\overline{G}^e$ both have
height $h$ over $\mathbf{F}_p$. The endomorphisms $F$, $V$ commute
with the action of $\mathbf{T}_\mathfrak{m}$, $F^{}$ acts on
$D(\overline{G}^e)$ by multiplication by the unit
$U_\p\cdot\langle\p\rangle^{-1}$ of $\mathbf{T}_\mathfrak{m}$, and
$V^{}$ acts by multiplication by the unit $U_\p$ of
$\mathbf{T}_\mathfrak{m}$.
\item The exact sequence $0\to G^0\to G\to G^e\to 0$ of
$p$-divisible groups over $F_\p^1$ gives a
$\Gal(\overline{F}_\p/F_\p)$-stable filtration $0\to T_pG^0\to T_p
G\to T_p G^e\to 0$. $\Gal(\overline{F}_\p/F_\p)$ acts on $T_pG^0$
by the character $\lambda(U_\p^{-1})\cdot\chi$, where $\chi$ is
the cyclotomic character and $\lambda(\cdot)$ is the unramified
character sending $\Frob_\p$ to $\cdot$, and on $T_pG^e$ by the
character
$\lambda(U_\p\cdot\langle\p\rangle^{-1})\cdot\chi^{2-k}$. Thus
there is a short exact sequence $0\to G^0[\mathfrak{m}]\to
G[\mathfrak{m}]\to G^e[\mathfrak{m}]\to 0$ over $\mathbf{Q}_p$,
with flat extensions to $\bigO_\p[\alpha]$, where
$G^e[\mathfrak{m}]$ does not necessarily denote the full
$\mathfrak{m}$-torsion in $G^e$, but rather the cokernel of the
map $G^0[\mathfrak{m}]\to G[\mathfrak{m}]$. The Galois group
$\Gal(\overline{F}_\p/F_\p)$ acts on $G^0[\mathfrak{m}]$ via the
character $\lambda(1/a_\p)\cdot\chi$ and on $G^e[\mathfrak{m}]$
via the character $\lambda(a_\p/\epsilon'(\p))\cdot\chi^{2-k}$.
\end{enumerate}
\end{thm}
\begin{proof}\enumerate \item The height of $G$ is equal to the
dimension of $V_pG=T_pG\otimes\mathbf{Q}_p$ as a
$\mathbf{Q}_p$-vector space. But by Lemma \ref{multone} $V_pG$ is
a free $\mathbf{T}_\mathfrak{m}\otimes\Q_p$-module of rank 2, so
the height $G$ is $2h$.

That $G$ is isomorphic to $G^\vee$ over $E_\mathfrak{np}$ follows
from the existence of a nondegenerate alternating form
$\langle,\rangle:T_pG\times T_pG\to T_p \mathbf{G}_m$ satisfying
$$\langle a^{\sigma_\l},b^{\sigma_\l}\rangle=\frac{\mathbf{N}\l}{\psi_{\pi'}(\l)}\langle
a,b\rangle$$ where $\sigma_\l$ is a Frobenius element at $\l$. The
existence of such a form can either be deduced as in \S 11 of
\cite{gro90} by modifying the Weil pairing on $J$ by an
Atkin-Lehner involution, or by (somewhat perversely) deducing its
existence from Theorem \ref{frob}. Then if $\l$ is trivial in the
ray class group mod $\mathfrak{np}$ we have $\langle
a^{\sigma_\l},b^{\sigma_\l}\rangle=\mathbf{N}\l\langle
a,b\rangle$; but such $\sigma_\l$ are dense in
$\Gal(\overline{E}/E_\mathfrak{np})$, so $G$ is isomorphic to
$G^\vee$ over $E_\mathfrak{np}$.

The proof that $G$ has good reduction over $F'_\p$ is very similar
to the argument in the proof of Proposition 12.9.1 in
\cite{gro90}. Let $A$ be the connected subgroup of points $P$ in
$J$ with $\sum_{a\in(\bigO_\p)^*}\langle a \rangle P=0$, and $B$
the connected subgroup of points fixed by the action of the group
$\{\langle a\rangle:a\in(\bigO_\p)^*\}$. The isogeny
\begin{align*}\phi:J &\to A\times B\\P&\mapsto\left((q-1)P-\sum_{a\in(\bigO_\p)^*}\langle a \rangle
P,\sum_{a\in(\bigO_\p)^*}\langle a \rangle P\right)\end{align*}
has degree prime to $p$, and thus induces an isomorphism on
$p$-divisible groups, because the composite with the natural
injection $A\times B\hookrightarrow J$ is just
$J\stackrel{q-1}{\To}J$. Because $k\neq 2$, the $p$-divisible
group of $G$ is a subgroup of the $p$-divisible group of $A$, so
it suffices to prove that $A$ has good reduction over $F'_\p$.
This follows from the arguments of \cite{dera} I.3.7 and V.3.2.

\item This follows from knowledge of the action of $U_\p$ on the components $\mathbf{M}'_{Ig,H'}$ and $\mathbf{M}_{Ig,H'}^{'(\sigma^{-1})}$ of the special fibre $\mathbf{M}'_{\operatorname{bal}.U_1(\p);\det=1,H'}$. As in Proposition 12.9.2 of \cite{gro90} we have
$$\overline{G}\cong\lim_\leftarrow\Jac(\mathbf{M}'_{Ig,H'})[\mathfrak{m}^n]\times\Jac(\mathbf{M}_{Ig,H'}^{'(\sigma^{-1})})[\mathfrak{m}^n].$$We
claim that in fact we have
$\overline{G}^m=\lim\limits_\longleftarrow\Jac(\mathbf{M}'_{Ig,H'})[\mathfrak{m}^n]$
and
$\overline{G}^e=\lim\limits_\longleftarrow\Jac(\mathbf{M}_{Ig,H'}^{'(\sigma^{-1})})[\mathfrak{m}^n]$.
Indeed, $U_\p$ is a unit in $\mathbf{T}_\mathfrak{m}$, and acts as
$\Ver^{}$ on $\mathbf{M}'_{Ig,H'}$ (note that this $\Ver$ is the
dual of the $p$-power Frobenius), whence
$\Jac(\mathbf{M}'_{Ig,H'})[\mathfrak{m}^n]$ is multiplicative (by
the standard theory of Dieudonn\'{e} modules). Similarly, one can
check that $U_\p$ acts as $\Frob^{}\cdot\langle\p\rangle$ on
$\mathbf{M}_{Ig,H'}^{'(\sigma^{-1})}$, so
$\Jac(\mathbf{M}_{Ig,H'}^{'(\sigma^{-1})})[\mathfrak{m}^n]$ is
\'{e}tale. That both subgroups have height $h$ follows from the
self-duality of $\overline{G}$ over $E_\mathfrak{np}$.
\item  The filtration $0\to T_pG^0\to T_p
G\to T_p G^e\to 0$ is stable under the action of
$\Gal(\overline{F}_\p/F'_\p)$, which acts (by the above) via the
characters $\lambda(U_\p^{-1})\cdot\chi$ on $T_pG^0$ and
$\lambda(U_\p)$ on $T_pG^e$. These characters are nonconjugate, so
the filtration is in fact stable under the action of
$\Gal(\overline{F}_\p/F_\p)$. To determine this action it suffices
to compute the action of $\Gal(F'_\p/F_\p)$; but this may be
accomplished just as in \cite{gro90}. It is easy to check that
$\Gal(F'_\p/F_\p)$ acts trivially on $\mathbf{M}'_{Ig,H'}$ and by
$\langle \alpha \rangle$ on $\mathbf{M}_{Ig,H'}^{'(\sigma^{-1})}$,
so by $\chi^{2-k}$ on $\overline{G}^e$, as required. The existence
of the short exact sequence $$0\to G^0[\mathfrak{m}]\to
G[\mathfrak{m}]\to G^e[\mathfrak{m}]\to 0$$follows at once, and
the action of $\Gal(\overline{F}_\p/F_\p)$ comes from the
observation that it acts semisimply by the results of \cite{blr}.

\end{proof}
\begin{corollary}Let $W$ be the two-dimensional vector space
underlying the representation $\overline{\rho}_f$. Then there is a
short exact sequence of $\Gal(\overline{F}_\p/F_\p)$-modules $0\to
V\to W\to V'\to 0.$ The group $\Gal(\overline{F}_\p/F_\p)$ acts on
$V$ by the character $\chi^{k-1}\lambda(\epsilon'(\p)/a_\p)$ and
on $V'$ by the unramified character $\lambda(a_\p)$. Equivalently,
there is a basis for $W$ such that $\Gal(\overline{F}_\p/F_\p)$
acts via the upper triangular matrices
$$\begin{pmatrix}\chi^{k-1}\lambda(\epsilon'(\p)/a_\p)&*\\0&\lambda(a_\p)\end{pmatrix}.$$
\end{corollary}
\begin{proof}This is immediate.\end{proof}

\begin{defn}Let $L$ be the field of definition of the
representation
$\overline{\rho}_f:\Gal(\overline{E}/E)\to\GL_2(\overline{\mathbf{F}}_p)$
i.e. the smallest field $L/\mathbf{F}_p$ such that
$\overline{\rho}_f$ factors through $GL_2(L)$.
\end{defn}
We have a realisation of
$\overline{\rho}_f\otimes(\epsilon'\chi^{k-2})^{-1}$ on
$G[\mathfrak{m}]$. The $L$-vector space scheme $G[\mathfrak{m}]$
sits in a short exact sequence over $F_\p$
\begin{equation}0\to G^0[\mathfrak{m}]\to G[\mathfrak{m}]\to
G^e[\mathfrak{m}] \to 0\label{exact}\tag{$\star$}\end{equation} of
$L$-vector space schemes, with $G^0[\mathfrak{m}]$ and
$G^e[\mathfrak{m}]$ both one-dimensional.

By their definitions, the $L$-vector space schemes in
\eqref{exact} all have canonical extensions to $\bigO_{F^1_\p}$.
In order to determine when $\overline{\rho}_f$ is tamely ramified,
we will determine when \eqref{exact} splits; in fact:

\begin{lemma}The following are equivalent:\begin{enumerate}\item The
sequence of $L$-vector space schemes in \eqref{exact} is uniquely
split over $F_\p$. \item The sequence of $L$-vector space schemes
over $\bigO_{F^1_\p}$ which extends \eqref{exact} is uniquely
split over $\bigO_{F^1_\p}$. \item The restriction of
$\overline{\rho}_f$ to $\Gal(\overline{F}_\p/F_\p)$ is
diagonalisable, and the sum of distinct characters
$\chi^{k-1}\lambda(\epsilon'(\p)/a_\p)$ and $\lambda(a_\p)$.
\end{enumerate}
\end{lemma}
\begin{proof}The equivalence of (1) and (3) is immediate. The
extension $F^1_\p/F_\p$ has degree prime to $p$, so (1) is
equivalent to the splitting of \eqref{exact} over $F^1_\p$, which
is obviously implied by (2). To establish that (1) implies (2), we
may base change to an \'etale extension $R$ of degree prime to $p$
of $\bigO_{F^1_\p}$, and check that a splitting of \eqref{exact}
over the quotient field $S$ of $R$ implies a splitting over $R$.

Choose $R$ so that $\lambda(1/a_\p)$,
$\lambda(a_\p/\epsilon'(\p))$ are trivial on
$\Gal(\overline{F}_\p/S)$. Then, as in the proof of Proposition
13.2 of \cite{gro90}, the $L$-vector space scheme
$G^e[\mathfrak{m}]$ is isomorphic to the \'etale vector space
scheme $L=L\otimes(\Z/p\Z)$ with trivial Galois action over $R$,
and $G^0[\mathfrak{m}]$ is isomorphic to the Cartier dual
$L^t=L^\vee\otimes\mu_p$ over $R$, where
$L^\vee=\operatorname{Hom}(L,\Z/p\Z)$. But we have a
Kummer-theoretic canonical isomorphism of $L$-vector spaces
$$\operatorname{Ext}_R(L,L^\vee\otimes\mu_p)\isoto
R^*/{R^*}^p\otimes L^\vee$$where $\operatorname{Ext}_R$ classifies
extensions in the category of $L$-vector space schemes.

Thus the sequence \eqref{exact} over $R$ gives a class in
$R^*/{R^*}^p\otimes L^\vee$ which is zero if and only if
\eqref{exact} splits; but $R^*/{R^*}^p$ injects into
$S^*/{S^*}^p\otimes L^\vee$, so a splitting over $S$ implies one
over $R$.
\end{proof}Now let $R$ denote the completion of the ring of integers in the maximal unramified extension of $\bigO_{F^1_\p}$. We now define a bilinear pairing $$q_f:(B^t)^e(L)\times
B^e(L)\to(R^*/{R^*}^p)\otimes L^\vee$$where
$B:=G[\mathfrak{m}]\otimes(\epsilon'\otimes\chi^{k-2})$, just as
in \S6 of \cite{cv92}. An element of $(B^t)^e(L)\times B^e(L)$
corresponds to a pair of homomorphisms
$\alpha:G^0[\mathfrak{m}]\to L^\vee\otimes\mu_p$ and $\beta:L\to
G^e[\mathfrak{m}]$, which give (via push-out and pull-back) a
homomorphism
$\alpha_*\beta^*:\operatorname{Ext}_R(G^e[\mathfrak{m}],G^0[\mathfrak{m}])\to\operatorname{Ext}_R(L,L^\vee\otimes\mu_p)\isoto
R^*/{R^*}^p\otimes L^\vee$, and the required element of
$R^*/{R^*}^p\otimes L^\vee$ is the image of the extension class of
$G[\mathfrak{m}]$ under $\alpha_*\beta^*$.

We have a map $\tr^\vee:L^\vee\to\mathbf{F}_p$ given by $h\mapsto
h(1)$, and a map $d\log:R^*\to\Omega^1_{R/\Z_p}$ given by
$a\mapsto da/a$, and thus a pairing $d\log
q_f:=(d\log\otimes\tr^\vee)\circ q_f:(B^t)^e(L)\times
B^e(L)\to\Omega^1_{R/\Z_p}$.
\begin{lemma}\label{328}If $\overline{\rho}_f$ is tamely ramified at $\p$,
then the pairing $d\log q_f$ is trivial.
\end{lemma}
\begin{proof}It suffices to show that $q_f$ is trivial; but as
noted above $q_f$ is trivial if and only if \eqref{exact} splits
if and only if $\overline{\rho}_f|_{\Gal(\overline{F}_\p/F_\p)}$
is diagonalisable if and only if $\overline{\rho}_f$ is tamely
ramified above $\p$.
\end{proof}
Since the $p$-divisible group $B$ is ordinary, we have (see \S
\ref{12}) the Serre-Tate pairing
$$q:T\overline{B}\times T\overline{B}^t\to 1+\pi R$$and thus a
pairing $$d\log q:T\overline{B}\times
T\overline{B}^t\to\Omega^1_{R/\bigO_\p^{nr}}$$which extends by
scalars to a pairing $$d\log
q:(T\overline{B}(\overline{\mathbf{F}}_p)\otimes_{\Z_p}R)\times(T\overline{B}^t(\overline{\mathbf{F}}_p)\otimes_{\Z_p}R)\to\Omega^1_{R/\bigO_\p^{nr}}.$$
Let $\overline{G}(j)$ denote the subgroup on which $(\bigO/\p)^*$
acts via $\langle\alpha\rangle=\alpha^j$. Then we have:
\begin{thm}\label{137}If $\alpha\in T\overline{B}(-k')\otimes R$, $\beta\in
T\overline{B}^t(k')\otimes R$, and
$\omega_\alpha|_{\mathbf{M}'_{Ig,H'}}=\omega_f$ (where
$\omega_\alpha=\alpha^*(dt/t)$), then $$d\log
q(\alpha,\beta)=\frac{u}{(k'-1)!a_\p}\left((w_\p^*\omega_\beta)|{\mathbf{M}'_{Ig,H'}},[M^{k'}\omega_f]\right)_{\mathbf{M}'_{Ig,H'}}\pi^{k'-1}d\pi+\dots$$where
$u$ is as in Theorem \ref{Upw}.\end{thm}
\begin{proof}This is very similar to Theorem 4.4 of
\cite{cv92}. From Theorem \ref{Upw} we have $[(\omega_f|U_\p
w_\p)|_{X'}/\pi^{k'}]=-\frac{\epsilon'(\p)}{k'!}[M^{k'}\omega_f]$,
and $\omega_\alpha|U_\p=a_\p\omega_\alpha\text{ mod }\pi$, so the
result follows from Theorem \ref{121}.\end{proof}

The Serre-Tate pairing $q:B[p](\overline{\mathbf{F}_p})\times
B^t[p](\overline{\mathbf{F}_p})\to R^*/{R^*}^p$ is related to
$q_f$ as follows:
\begin{lemma}Suppose that $\alpha\in
B[p](\overline{\mathbf{F}}_p)$ and $\beta\in
B(\overline{\mathbf{F}}_p)$. Then
$$(1\otimes\tr^\vee)q_f(\alpha\text{ mod
}\mathfrak{m},\beta)=q(\alpha,\beta)\text{ mod }{R^*}^p.$$
\end{lemma}
\begin{proof}This follows from the definitions of $q_f$, $q$.
\end{proof}
We have a map $T(B^e)\to\operatorname{H}^0({\mathbf{M}'}^\text{can}_{\operatorname{bal}.U_1({\p}),H'},\Omega^1_{{\mathbf{M}'}^\text{can}_{\operatorname{bal}.U_1({\p}),H'}/R})
\to\operatorname{H}^0(\mathbf{M}'_{Ig,H'},\Omega^1_{\mathbf{M}'_{Ig,H'}/\overline{\mathbf{F}}_p})$. Let
$\beta_f$ be the element of
$B^e(\overline{\Q}_p)\otimes_{\mathbf{F}_p}\overline{\mathbf{F}}_p$
corresponding to $\omega_f$ via this map. Then
\begin{prop}\label{331}$[M^{k'}\omega_f]=0$ if and only if $d\log
q_f(\alpha,\beta_f)=0$ for all
$\alpha\in(B^t)^e(\overline{\Q}_p)\otimes_{\mathbf{F}_p}\overline{\mathbf{F}}_p$.\end{prop}
\begin{proof}By Theorem \ref{137} we have $d\log q_f(\alpha,\beta_f)=0$
for all
$\alpha\in(B^t)^e(\overline{\Q}_p)\otimes_{\mathbf{F}_p}\overline{\mathbf{F}}_p$
if and only if $(\eta,[M^{k'}\omega_f])_{\mathbf{M}'_{Ig,H'}}=0$ for
all $\eta\in\operatorname{H}^0(\mathbf{M}'_{Ig,H'}
,\Omega^1_{\mathbf{M}'_{Ig,H'}/\overline{\mathbf{F}}_p})$. But by
Theorem \ref{unitroot} $[M^{k'}\omega_f]$ is in the unit root
subspace of
$\operatorname{H}^1_{dR}(\mathbf{M}'_{Ig,H'}/\overline{\mathbf{F}}_p)$,
which has trivial intersection with the space of global
differentials, which is a maximal isotropic subspace for the pairing
$(,)_{\mathbf{M}'_{Ig,H'}}$. So if
$(\eta,[M^{k'}\omega_f])_{\mathbf{M}'_{Ig,H'}}=0$ for
$\eta\in\operatorname{H}^0(\mathbf{M}'_{Ig,H'},\Omega^1_{\mathbf{M}'_{Ig,H'}/\overline{\mathbf{F}}_p})$
then $[M^{k'}\omega_f]=0$. The converse is clear.\end{proof}

\section{The Main Theorem}\label{final}\subsection{Proof of the main result} From Theorem \ref{313}, Lemma \ref{328} and Proposition \ref{331} we see
that we have constructed a companion form $\omega_{f'}$ as a
differential. Note that $\omega_{f'}$ has character (of
$(\bigO_\p/\p)^\times)$) $\omega_\p^{-k}$ at $\p$, and
$\overline{\rho}_{\pi'}\cong\overline{\rho}_{\pi}$. Furthermore, it
has character $\omega_{\p_i}^{-k'}$ at $\p_i$ for $\p_i\neq\p$,
because the diamond operators at $\p_i$ act trivially on $a$ and
thus on $\omega_A$. We now reverse the process that associated
$\omega_f$ to $\pi$, associating a mod $p$ Hilbert modular form
$\pi'$ to $\omega_{f'}$. Firstly, by the Deligne-Serre lemma there
is a characteristic zero differential $\omega_{F'}$ whose Hecke
eigenvalues lift those of $\omega_{f'}$. This differential
corresponds to an automorphic form $\pi_2$ on $G'(\mathbf{A}_\Q)$
such that $\psi'_\infty$ is of weight 2, where
$BC(\pi_2)=(\psi',\Pi')$. Put $\pi'_E=JL(\Pi')$. Note that
$\pi_E'\cong(\pi_E')^\vee\circ c$; if now $\eta'$ is a character
such that ${\eta'}^c/\eta'=\chi_{\pi_E'}$, we have
\begin{align*}(\pi_E'\otimes\eta')\circ
c&\cong(\pi_E')^c\otimes{\eta'}^c\\&\cong(\pi_E')^\vee\otimes{\eta'}^c\\&\cong\pi_E'\otimes\chi_{\pi_E'}^{-1}\otimes{\eta'}^c\\&\cong
\pi_E'\otimes\eta',\end{align*} so there exists $\pi'$ on
$\GL_2(\mathbf{A}_F)$ such that $BC_{E/F}(\pi')=\pi_E'\otimes\eta'$.
Again, $(\eta')_\infty$ is trivial, and we choose $\eta'$ so that
$\eta'_\p$ is trivial. We wish to check that $\pi'$ is ordinary at
all primes dividing $p$; this follows from

\begin{thm}Let $\Pi$ be a weight two Hilbert modular form of level
$\n p$ and character  $\epsilon$ (a strict ray class character of
conductor $\n p$), and suppose it has character $\omega_{\p_i}^{j}$
at $\p_i$, with $0<j<p-2$. Let the slope of $\Pi$ be the $p$-adic
valuation of the eigenvalue $a_{\p_i}$ of $U_{\p_i}$ on $\Pi$. Then:
\begin{itemize}\item If $\Pi$ has slope 0, then
$\overline{\rho}_{\Pi}|_{I_{\p_i}}\cong\left(%
\begin{array}{cc}
  \omega_{\p_i}^{j+1} & * \\
  0 & 1 \\
\end{array}%
\right);$
\item If $\Pi$ has slope
1, then
$\overline{\rho}_{\Pi}|_{I_{\p_i}}\cong\left(%
\begin{array}{cc}
  \omega_{\p_i} & * \\
  0 & \omega_{\p_i}^j \\
\end{array}%
\right);$
\item If $\Pi$ has slope
in the interval(0,1), then $\overline{\rho}_{\Pi}|_{I_{\p_i}}\cong
\omega_2^{1+j}\oplus\omega_2^{p(1+j)}$, where $\omega_2$ is a
fundamental character of niveau 2 associated to $\p_i$.
\end{itemize}
\end{thm}
\begin{proof}The proof is very similar to that of Proposition 6.17
of \cite{sav04}. In fact, the proof, being purely local (and recall
that we are in the case where $p$ is totally split in $F$), is
identical once one has a generalisation of the results of
\cite{sai97} to the Hilbert case, for which see Theorem 1 of
\cite{sai03}
\end{proof}
We now repeat the above arguments at all other primes of $F$
dividing $p$, until we obtain a weight 2 level $\n p$ Hilbert
modular form $\pi'$ with character $\omega_{\p_i}^{-k'}$ at all
$\p_i|p$, with
$\overline{\rho}_{\pi'}|_{\Gal(\overline{E}/E)}\cong\overline{\rho}_\pi|_{\Gal(\overline{E}/E)}$.

This does not, of course, guarantee that we have
$\overline{\rho}_{\pi'}\cong\overline{\rho}_\pi$. However:
\begin{thm}\label{penultimate}There exists a weight $2$ level $\mathfrak{n}p$ Hilbert
modular form $\pi'$ with character $\omega_{\p_i}^{-k'}$ at all
$\p_i|p$ such that the representation
$\overline{\rho}_{\pi'}:\Gal(\overline{F}/F)\to\GL_2(\overline{\mathbf{F}}_p)$
satisfies $\overline{\rho}_{\pi'}\cong\overline{\rho}_\pi$.
\end{thm}
\begin{proof}Suppose not. For any imaginary quadratic extension $K/\Q$ as above, and CM extension
$E=FK$ of $F$, we can construct a $\pi'_K$ as above, satisfying
$\overline{\rho}_{\pi'_K}|_{\Gal(\overline{E}/E)}\cong\overline{\rho}_f|_{\Gal(\overline{E}/E)}$.
Since there are only finitely many mod $p$ weight $2$ level
$\mathfrak{n}p$ Hilbert modular forms, there are only a finite
number of ``candidate'' automorphic forms $\pi_i$. For each
$\pi_i$ there must (by the Cebotarev density theorem and our
assumption that there is no such $\pi'$ with
$\overline{\rho}_{\pi'}\cong\overline{\rho}_\pi$) be infinitely
many places at which $\pi_i$ is unramified principal series with
the ``wrong'' characters mod $p$. In particular, we may choose
such a place $v_i$ for each $\pi_i$; then choosing $K/\Q$ to split
at the places of $\Q$ lying below all of the $v_i$ gives a
contradiction.
\end{proof}

\begin{thmn}[Theorem A]Let $F$ be a totally real field in which an odd prime $p$
splits completely. Let $\pi$ be a mod $p$ Hilbert modular form of
parallel weight $2<k<p$ and level $\n$, with $\n$ coprime to $p$.
Suppose that $\pi$ is ordinary at all primes $\p|p$, and that the
mod $p$ representation
$\overline{\rho}_\pi:\Gal(\overline{F}/F)\to\GL_2(\overline{\F}_p)$
is irreducible and is tamely ramified at all primes $\p|p$. Then
there is a companion form $\pi'$ of parallel weight $k'=p+1-k$ and
level $\n$ satisfying
$\overline{\rho}_{\pi'}\cong\overline{\rho}_\pi\otimes\chi^{k'-1}$.\end{thmn}
\begin{proof}Firstly, we deal with the case where $[F(\zeta_p):F]=2$ and $\overline{\rho}_\pi|_{\Gal(\overline{F}/F(\zeta_p))}$
is reducible. In this case we can directly construct a companion
form, in a similar fashion to \cite{wie04}. Since
$\overline{\rho}_\pi|_{\Gal(\overline{F}/F(\zeta_p))}$ is
reducible, $\overline{\rho}$ is induced from a character $\psi$ on
$\Gal(\overline{F}/F(\zeta_p))$. As in Lemma 2 of \cite{wie04} we
take the Teichmuller lift
$\widetilde{\psi}:\Gal(\overline{F}/F(\zeta_p))\to\bigO_{F[\zeta_p]}^\times$,
so that
$\widetilde{\rho}:=\operatorname{Ind}_{\Gal(\overline{F}/F(\zeta_p))}^{\Gal(\overline{F}/F)}\widetilde{\psi}$
is an odd (as $p>2$) lift of $\overline{\rho}$, corresponding to a
modular form of weight 1. By hypothesis
$\overline{\rho}|_{\Gal(\overline{F}_\p/F_\p)}$ is the direct sum
of two characters, one of which is unramified; but then it is easy
to see that ${\widetilde{\psi}}^2$ must be unramified, which
yields $k-1=p-1$ or $k-1=(p-1)/2$. In the first case, the weight
one form is unramified principal series at $p$, and is thus the
required companion form; and in the second case we twist with the
quadratic character whose restriction to
$\Gal(\overline{F}/F(\zeta_p))$ is trivial, and then use Hida
theory to move to the required companion form in weight $(p+1)/2$.

In the general case, we must prove that the mod $p$, level
$\mathfrak{n}p$ weight 2 modular form $\pi'$ constructed in Theorem
\ref{penultimate} is also of weight $k'$ and level $\mathfrak{n}$.
But again, it follows from the Hida theory in \cite{wil88} that we
have an ordinary form of weight $k'$ and level $\mathfrak{n}p$, and
it remains to check that this form cannot be new at any $\p_i|p$.
This is, however, easy; as the character at $\p_i$ is unramified
(because we are in weight $k'$), $\pi'_{\p_i}$ would have to be
special, which immediately contradicts the fact that it is ordinary,
provided that $k'\neq 2$; but in this case we may use Theorem 6.2 of
\cite{jar04}.
\end{proof}

\def\cprime{$'$}
\providecommand{\bysame}{\leavevmode\hbox
to3em{\hrulefill}\thinspace}
\providecommand{\MR}{\relax\ifhmode\unskip\space\fi MR }

\providecommand{\MRhref}[2]{
  \href{http://www.ams.org/mathscinet-getitem?mr=#1}{#2}
} \providecommand{\href}[2]{#2}


\begin{thebibliography}{Car86b}

\bibitem[BLR91]{blr}
Nigel Boston, Hendrik~W. Lenstra, Jr., and Kenneth~A. Ribet,
\emph{Quotients of
  group rings arising from two-dimensional representations}, C. R. Acad. Sci.
  Paris S\'er. I Math. \textbf{312} (1991), no.~4, 323--328.

\bibitem[BR89]{br89}
D.~Blasius and J.~Rogawski, \emph{Galois representations for
{H}ilbert modular
  forms}, Bull. Amer. Math. Soc. (N.S.) \textbf{21} (1989), no.~1, 65--69.

\bibitem[BT99]{bt}
Kevin Buzzard and Richard Taylor, \emph{Companion forms and weight
one forms},
  Ann. of Math. (2) \textbf{149} (1999), no.~3, 905--919.

\bibitem[Car86a]{car861}
Henri Carayol, \emph{Sur la mauvaise r\'eduction des courbes de
{S}himura},
  Compositio Math. \textbf{59} (1986), no.~2, 151--230.

\bibitem[Car86b]{car862}
\bysame, \emph{Sur les repr\'esentations {$l$}-adiques associ\'ees
aux formes
  modulaires de {H}ilbert}, Ann. Sci. \'Ecole Norm. Sup. (4) \textbf{19}
  (1986), no.~3, 409--468.

\bibitem[CL99]{cl99}
Laurent Clozel and Jean-Pierre Labesse, \emph{Changement de base
pour les
  repr\'esentations cohomologiques de certains groupes unitaires}, Ast\'erisque
  (1999), no.~257.

\bibitem[CV92]{cv92}
Robert~F. Coleman and Jos{\'e}~Felipe Voloch, \emph{Companion forms
and
  {K}odaira-{S}pencer theory}, Invent. Math. \textbf{110} (1992), no.~2,
  263--281.

\bibitem[DR73]{dera}
P.~Deligne and M.~Rapoport, \emph{Les sch\'emas de modules de
courbes
  elliptiques}, Modular functions of one variable, II (Proc. Internat. Summer
  School, Univ. Antwerp, Antwerp, 1972), Springer, Berlin, 1973, pp.~143--316.
  Lecture Notes in Math., Vol. 349.

\bibitem[Dri76]{dri76}
V.~G. Drinfel{\cprime}d, \emph{Elliptic modules}, Math. USSR-Sb.
(N.S.)
  \textbf{23(4)} (1976), 561--592.

\bibitem[DS74]{delser}
Pierre Deligne and Jean-Pierre Serre, \emph{Formes modulaires de
poids {$1$}},
  Ann. Sci. \'Ecole Norm. Sup. (4) \textbf{7} (1974), 507--530 (1975).

\bibitem[FC90]{cf90}
Gerd Faltings and Ching-Li Chai, \emph{Degeneration of abelian
varieties},
  Ergebnisse der Mathematik und ihrer Grenzgebiete (3) [Results in Mathematics
  and Related Areas (3)], vol.~22, Springer-Verlag, Berlin, 1990, With an
  appendix by David Mumford.

\bibitem[Fuj99]{fuj99}
Kazuhiro Fujiwara, \emph{Level optimisation in the totally real
case}, 1999.

\bibitem[Gee]{gee06}
Toby Gee, \emph{Companion forms over totally real fields II}, to appear in Duke Mathematical Journal.

\bibitem[Gro67]{egaiv}
A.~Grothendieck, \emph{\'{E}l\'ements de g\'eom\'etrie alg\'ebrique.
{IV}.
  \'{E}tude locale des sch\'emas et des morphismes de sch\'emas {IV}}, Inst.
  Hautes \'Etudes Sci. Publ. Math. (1967), no.~32, 361.

\bibitem[Gro90]{gro90}
Benedict~H. Gross, \emph{A tameness criterion for {G}alois
representations
  associated to modular forms (mod {$p$})}, Duke Math. J. \textbf{61} (1990),
  no.~2, 445--517.

\bibitem[Har00]{har00}
Michael Harris, \emph{The {L}ocal {L}anglands {C}orrespondence:
Notes of
  ({H}alf) a {C}ourse at the {IHP} {S}pring 2000}, 2000.

\bibitem[HT01]{ht01}
Michael Harris and Richard Taylor, \emph{The geometry and cohomology
of some
  simple {S}himura varieties}, Annals of Mathematics Studies, vol. 151,
  Princeton University Press, Princeton, NJ, 2001, With an appendix by Vladimir
  G. Berkovich.

\bibitem[HT02]{ht02}
\bysame, \emph{Regular models of certain {S}himura varieties}, Asian
J. Math.
  \textbf{6} (2002), no.~1, 61--94.

\bibitem[Ill85]{ill85}
Luc Illusie, \emph{D\'eformations de groupes de {B}arsotti-{T}ate
(d'apr\`es
  {A}. {G}rothendieck)}, Ast\'erisque (1985), no.~127, 151--198, Seminar on
  arithmetic bundles: the Mordell conjecture (Paris, 1983/84).

\bibitem[Jar99]{jar99}
Frazer Jarvis, \emph{Mazur's principle for totally real fields of
odd degree},
  Compositio Math. \textbf{116} (1999), no.~1, 39--79.

\bibitem[Jar04]{jar04}
\bysame, \emph{Correspondences on {S}himura curves and {M}azur's
{P}rinciple
  above $p$}, Pacific J. Math. \textbf{213} (2004), no.~2, 267--280.

\bibitem[JL70]{jaclan}
H.~Jacquet and R.~P. Langlands, \emph{Automorphic forms on {${\rm
GL}(2)$}},
  Springer-Verlag, Berlin, 1970, Lecture Notes in Mathematics, Vol. 114.

\bibitem[Kas04]{kas04}
Payman~L. Kassaei, \emph{{$P$}-adic modular forms over {S}himura
curves over
  totally real fields}, Compos. Math. \textbf{140} (2004), no.~2, 359--395.

\bibitem[Kat89]{kat89}
Kazuya Kato, \emph{Logarithmic structures of {F}ontaine-{I}llusie},
Algebraic
  analysis, geometry, and number theory (Baltimore, MD, 1988), Johns Hopkins
  Univ. Press, Baltimore, MD, 1989, pp.~191--224.

\bibitem[Kat81]{kat81}
N.~Katz, \emph{Serre-{T}ate local moduli}, Algebraic surfaces
(Orsay,
  1976--78), Lecture Notes in Math., vol. 868, Springer, Berlin, 1981,
  pp.~138--202.



\bibitem[KM85]{km}
Nicholas~M. Katz and Barry Mazur, \emph{Arithmetic moduli of
elliptic curves},
  Annals of Mathematics Studies, vol. 108, Princeton University Press,
  Princeton, NJ, 1985.

\bibitem[Kot92]{ko92}
Robert~E. Kottwitz, \emph{On the {$\lambda$}-adic representations
associated to
  some simple {S}himura varieties}, Invent. Math. \textbf{108} (1992), no.~3,
  653--665.

\bibitem[Lan80]{lan80}
Robert~P. Langlands, \emph{Base change for {${\rm GL}(2)$}}, Annals
of
  Mathematics Studies, vol.~96, Princeton University Press, Princeton, N.J.,
  1980.

\bibitem[LR04]{lr04}
Joshua Lansky and A.~Raghuram, \emph{Conductors and {N}ewforms for
${U}(1,1)$},
  Proceedings of the Indian Academy of Sciences. \textbf{114}, No. 4, 319--343, (2004).  

\bibitem[MW84]{mw84}
B.~Mazur and A.~Wiles, \emph{Class fields of abelian extensions of
{${\bf
  Q}$}}, Invent. Math. \textbf{76} (1984), no.~2, 179--330.

\bibitem[Sai97]{sai97}
Takeshi Saito, \emph{Modular forms and {$p$}-adic {H}odge theory},
Invent.
  Math. \textbf{129} (1997), no.~3, 607--620.

\bibitem[Sai03]{sai03}
\bysame, \emph{Hilbert modular forms and p-adic {H}odge theory},
2003.

\bibitem[Sav04]{sav04}
David Savitt, \emph{On a conjecture of {C}onrad, {D}iamond and
{T}aylor}, Duke Mathematical Journal \textbf{128} (2005), no. 1, 141-197.

\bibitem[Tay89]{tay89}
Richard Taylor, \emph{On {G}alois representations associated to
{H}ilbert
  modular forms}, Invent. Math. \textbf{98} (1989), no.~2, 265--280.

\bibitem[Wie04]{wie04}
Gabor Wiese, \emph{Dihedral {G}alois representations and {K}atz
modular forms}, Documenta Mathematica \textbf{9} (2004), 123-133.

\bibitem[Wil80]{wil80}
Andrew Wiles, \emph{Modular curves and the class group of {${\bf
Q}(\zeta
  \sb{p})$}}, Invent. Math. \textbf{58} (1980), no.~1, 1--35.

\bibitem[Wil88]{wil88}
A.~Wiles, \emph{On ordinary {$\lambda$}-adic representations
associated to
  modular forms}, Invent. Math. \textbf{94} (1988), no.~3, 529--573.

\end{thebibliography}
\end{document}